\documentclass[11pt,a4paper]{article}

\usepackage[plainpages = false]{hyperref}
\usepackage{amsmath}
\usepackage{amsthm}
\usepackage{thmtools, thm-restate}
\usepackage{amssymb}
\usepackage[capitalize]{cleveref}
\usepackage[british]{babel}
\usepackage{tikz}
\usepackage{tikz-cd}
\usetikzlibrary{babel}
\usepackage{csquotes}
\usepackage{enumerate}
\usepackage{tkz-berge}
\usepackage{scalerel}
\usepackage{graphicx,xcolor,textpos}
\usepackage{url}
\usepackage{xstring}
\usepackage{ifthen}
\usepackage[backend = biber, style = numeric, giveninits]{biblatex}
\usepackage{a4wide}
\usepackage{subcaption}

\renewbibmacro*{issue+date}{%
  \setunit*{\addcomma\space}
    \iffieldundef{issue}
      {\usebibmacro{date}}
      {\printfield{issue}%
       \setunit*{\addspace}%
       \usebibmacro{date}}
  \newunit}

\renewbibmacro*{journal+issuetitle}{%
  \usebibmacro{journal}%
  \setunit*{\addspace}%
  \iffieldundef{series}
    {}
    {\newunit
     \printfield{series}%
     \setunit{\addspace}}%
  \setunit*{\addcomma\space}
  \usebibmacro{volume+number+eid}%
  \newunit}

\DeclareBibliographyDriver{article}{%
  \usebibmacro{bibindex}%
  \usebibmacro{begentry}%
  \usebibmacro{author/translator+others}%
  \setunit{\labelnamepunct}\newblock
  \usebibmacro{title}%
  \newunit
  \printlist{language}%
  \newunit\newblock
  \usebibmacro{byauthor}%
  \newunit\newblock
  \usebibmacro{bytranslator+others}%
  \newunit\newblock
  \printfield{version}%
  \newunit\newblock
  \usebibmacro{journal+issuetitle}%
  \newunit
  \usebibmacro{byeditor+others}%
  \newunit
  \usebibmacro{note+pages}%
  \setunit{\addspace}
  \usebibmacro{issue+date}
  \setunit{\addcolon\space}
  \usebibmacro{issue}
  \newunit\newblock
  \iftoggle{bbx:isbn}
    {\printfield{issn}}
    {}%
  \newunit\newblock
  \usebibmacro{doi+eprint+url}%
  \newunit\newblock
  \usebibmacro{addendum+pubstate}%
  \setunit{\bibpagerefpunct}\newblock
  \usebibmacro{pageref}%
  \usebibmacro{finentry}}

\renewbibmacro{in:}{}

\renewrobustcmd*{\bibinitdelim}{\,}

\DeclareNameAlias{default}{last-first}

\DeclareFieldFormat[article,periodical]{number}{\mkbibparens{#1}}

\renewbibmacro*{volume+number+eid}{%
  \printfield{volume}%
  \printfield{number}%
  \setunit{\addcomma\space}%
  \printfield{eid}}

\DeclareFieldFormat[article]{pages}{#1}

\newcommand{\ceil}[1]{\left\lceil{#1}\right\rceil}
\newcommand{\cl}[1]{\overline{#1}}
\newcommand{\conj}[1]{\sim_{#1}}

\newcommand{\eg}{e.g.\ }

\newcommand{\G}{\mathcal{G}}
\newcommand{\gen}[1]{\langle #1 \rangle}
\newcommand{\grp}[2]{\langle #1 \mid #2 \rangle}
\newcommand{\ie}{i.e.\ }

\newcommand{\inv}[1]{{#1}^{-1}}

\newcommand{\liebarg}[2]{[#1, #2]}
\newcommand{\lieb}{\liebarg{\cdot}{\cdot}}

\newcommand{\N}{\mathbb{N}}
\newcommand{\normcl}[1]{\langle\langle #1 \rangle \rangle}

\newcommand{\pres}[2]{\langle #1 \mid #2 \rangle}
\newcommand{\Reid}{\mathcal{R}}

\newcommand{\Rinf}{R_\infty}

\newcommand{\Vchar}{\Vtext{char}}
\newcommand{\Vmax}{\Vtext{max}}
\newcommand{\Vnottrans}{\Vtext{\overline{\tau}}}
\newcommand{\Vtext}[1]{V_{\mathrm{#1}}}
\newcommand{\Vtrans}{\Vtext{\tau}}
\newcommand{\wrt}{w.r.t.\ }
\newcommand{\Z}{\mathbb{Z}}

\renewcommand{\phi}{\varphi}

\let\originalleft\left
\let\originalright\right
\renewcommand{\left}{\mathopen{}\mathclose\bgroup\originalleft}
\renewcommand{\right}{\aftergroup\egroup\originalright}

\DeclareMathOperator*{\Ast}{\scalerel*{\ast}{\textstyle\sum}}
\DeclareMathOperator{\Aut}{Aut}
\DeclareMathOperator{\Autnottrans}{\Aut_{\overline{\tau}}}
\DeclareMathOperator{\Autip}{\Aut_{\iota,p}}

\DeclareMathOperator{\GL}{GL}

\DeclareMathOperator{\Id}{Id}

\DeclareMathOperator{\lcm}{lcm}

\DeclareMathOperator{\len}{len}
\DeclareMathOperator{\SL}{SL}
\DeclareMathOperator{\Span}{span}
\DeclareMathOperator{\Spec}{Spec}

\DeclareMathOperator*{\Times}{\scalerel*{\times}{\textstyle\sum}}

\declaretheorem[style=definition, name = Definition, numberwithin=subsection]{defin}
\declaretheorem[style=definition, name = Example, sibling=defin]{example}
\declaretheorem[name = Theorem, sibling=defin]{theorem}
\declaretheorem[name = Lemma, sibling=defin]{lemma}
\declaretheorem[name = Proposition, sibling=defin]{prop}
\declaretheorem[name = Corollary, sibling=defin]{cor}
\declaretheorem[style=remark, name = Remark, numbered = no]{remark}
\declaretheorem[name = Conjecture, numbered = no]{conjecture}
\declaretheorem[name = Theorem, numbered = no]{theorem*}

\numberwithin{equation}{section}

\crefname{prop}{Proposition}{Propositions}
\crefname{cor}{Corollary}{Corollaries}

\begin{document}
\begingroup
\begin{center}
	\LARGE{The \(\Rinf\)-property for right-angled Artin groups}	\\[.5em]
	\Large{Karel Dekimpe\footnote{Research supported by long term structural funding -- Methusalem grant of the Flemish Government.} and Pieter Senden\footnote{Researcher funded by FWO PhD-fellowship fundamental research (file number: 1112520N).}}
\end{center}
\endgroup

\begin{abstract}
Given a group \(G\) and an automorphism \(\phi\) of \(G\), two elements \(x, y \in G\) are said to be \(\phi\)-conjugate if \(x = g y \inv{\phi(g)}\) for some \(g \in G\). The number of equivalence classes is the Reidemeister number \(R(\phi)\) of \(\phi\), and if \(R(\phi) = \infty\) for all automorphisms of \(G\), then \(G\) is said to have the \(\Rinf\)-property. 

A finite simple graph \(\Gamma\) gives rise to the right-angled Artin group \(A_{\Gamma}\), which has as generators the vertices of \(\Gamma\) and as relations \(vw = wv\) if and only if \(v\) and \(w\) are joined by an edge in \(\Gamma\). We conjecture that all non-abelian right-angled Artin groups have the \(\Rinf\)-property and prove this conjecture for several subclasses of right-angled Artin groups.
\end{abstract}
\let\thefootnote\relax
\let\thefootnote\relax
\section{Twisted conjugacy and Reidemeister numbers}	\label{sec:twistedConjugacy}
\addtocounter{subsection}{1}
	Let $G$ be a group and $\phi: G \to G$ be an automorphism. For $x, y \in G$, we say that $x$ and $y$ are \emph{$\phi$-conjugate} and write $x \conj{\phi} y$ if there exists a $g \in G$ such that $x = g y \inv{\phi(g)}$. The equivalence class of $x$ is denoted by $[x]$, or $[x]_\phi$ for clarity if there are multiple automorphisms involved.

	We define $\Reid[\phi]$ to be the set of all $\conj{\phi}$-equivalence classes and the \emph{{Reidemeister} number $R(\phi)$ of $\phi$} as the cardinality of $\Reid[\phi]$. Note that $R(\phi) \in \N_0 \cup \{\infty\}$. Finally, we define the \emph{Reidemeister spectrum} to be
	\(
		\Spec_R(G) := \{ R(\phi) \mid \phi \in \Aut(G)\}.
	\)
	We say that \emph{$G$ has the $\Rinf$-property}, also denoted as $G \in \Rinf$, if $\Spec_R(G) = \{\infty\}$. We say that $G$ has \emph{full Reidemeister spectrum} if $\Spec_R(G) = \N_0 \cup \{\infty\}$.
	
The notion of Reidemeister number arises from Nielsen fixed-point theory, where its topological analog serves as a count of the number of fixed point classes of a continuous self-map, and is strongly related to the algebraic one defined above, see \cite{Jiang}.

It has been proven for several (classes of) groups that they possess the \(\Rinf\)-property, \eg Baumslag-Solitar groups \cite{BaumslagSolitar} and their generalisations \cite{Levitt}, extensions of \(\SL(n, \Z)\) and \( \GL(n, \Z)\) by a countable abelian group \cite{MubeenaSankaran}, and Thompson's group \cite{ThompsonGroup}. We refer the reader to \cite{FelshtynNasybullov} for a more exhaustive list of groups having the  \(\Rinf\)-property.  \newline

In this article, we study the Reidemeister spectrum of right-angled Artin groups, RAAGs for short. Given a graph \(\Gamma\) with the set of vertices \(V\), the RAAG associated to it is the group
\[
	A_{\Gamma} = \grp{V}{[v, w] \text{ if $v, w \in V$ are joined by an edge in $\Gamma$}}.
\]
Extreme cases of RAAGs include free groups and free abelian groups, coming from edgeless and complete graphs, respectively. From \cite[Theorem 3]{Felshtyn} (see also \cite{DekimpeGoncalves}), it readily follows that all non-abelian free groups of finite rank have the \(\Rinf\)-property. On the other hand, \(\Spec_{R}(\Z) = \{2, \infty\}\) and \(\Spec_{R}(\Z^{n})=\N_0\cup \{\infty\}\) for \(n \geq 2\) (see \eg \cite{Romankov}). For groups closely related to right-angled Artin groups, several results regarding the \(\Rinf\)-property have been obtained, \eg A.\ Juh\'asz has proven that certain Artin groups which are not right-angled Artin groups possess the \(\Rinf\)-property \cite{Juhasz} and recently, T.\ K.\ Naik, N.\ Nanda and M.\ Singh have showed that twin groups, a subfamily of the right-angled Coxeter groups, all possess the \(\Rinf\)-property \cite{NaikNandaSingh}.

We suspect that, amongst all RAAGs, only the free abelian ones do not possess the \(\Rinf\)-property:
\begin{conjecture}
	Let $\Gamma(V, E)$ be a finite non-complete graph, \ie $V$ is finite and there are two (distinct) vertices not joined by an edge. Then \(A_{\Gamma}\) has the \(\Rinf\)-property.
\end{conjecture}

We first reduce the conjecture to graphs belonging to three specific classes, after which we prove the conjecture for one of these classes and for several subclasses of the other two. 

We start by recalling two ways of proving that an automorphism has infinite Reidemeister number.

\begin{defin}
	Let $G, H$ be two groups. If $H \cong G / N$ for some characteristic subgroup $N$ of $G$, we call $H$ a \emph{characteristic quotient of $G$}.
\end{defin}
The following result is well-known, see \eg \cite[Lemma 2.1]{MubeenaSankaran}.

\begin{lemma}	\label{cor:characteristicQuotientRinf2}
	Let $G, H$ be two groups. If $H$ is a characteristic quotient of $G$ and $H \in \Rinf$, then $G \in \Rinf$.
\end{lemma}
Reidemeister numbers also behave nicely under conjugation. For elements \(a, b\) of a group \(G\), we put \(a^{b} := \inv{b} a b\).
\begin{prop}	\label{prop:conjugateEndomorphisms}
	Let $G$ be a group and $\phi, \psi \in \Aut(G)$. Then $R(\phi) = R(\phi^\psi)$.
\end{prop}
\begin{proof}

	Define $\rho: \Reid[\phi] \to \Reid[\phi^\psi]: [g]_\phi \to [\inv{\psi}(g)]_{\phi^\psi}$. As $\inv{\psi}$ is surjective, so is $\rho$, and
	\begin{align*}
		x \conj{\phi} y 	&\iff	\exists g \in G: x = gy\inv{\phi(g)}	\\
					&\iff	\exists g \in G: \inv{\psi}(x) = \inv{\psi}(g) \inv{\psi}(y) \inv{\inv{\psi}(\phi(\psi(\inv{\psi}(g))))}	\\
					&\iff \exists g \in G: \inv{\psi}(x) = \inv{\psi}(g) \inv{\psi}(y) \inv{\phi^\psi(\inv{\psi}(g))}	\\
					&\iff \inv{\psi}(x) \conj{\phi^\psi} \inv{\psi}(y),
	\end{align*}
	which shows well-definedness and injectivity of $\rho$.
\end{proof}

Next, we recall the definition and some properties of the lower central series of a group.
\begin{defin}

	Let \(G\) be a group. The \emph{lower central series} of $G$ is defined as follows: put $\gamma_1(G) = G$ and inductively define \(\gamma_{i + 1}(G) = [\gamma_i(G), G]\) for $i \geq 1$.
\end{defin}

Each term in the lower central series is a characteristic subgroup, the quotients \(\gamma_i(G) / \gamma_{i + 1}(G)\) are all abelian and we will refer to them as the \emph{factors of the lower central series}. We can put these together to form the associated Lie ring of \(G\):

\begin{defin}
	Let $G$ be a group. The \emph{Lie ring associated to $G$} is the abelian group
	\[
		L(G) := \bigoplus_{i = 1}^\infty L_i(G), \quad \text{ where }	L_i(G) = \frac{\gamma_i(G)}{\gamma_{i + 1}(G)},
	\]
	and $L(G)$ is equipped with the following Lie bracket $\lieb{}_{L}$: for $g \in \gamma_i(G)$ and $h \in \gamma_j(G)$, we define
	\[
		[g\gamma_{i + 1}(G), h\gamma_{j +1}(G)]_L := [g, h]\gamma_{i + j + 1}(G),
	\]
	and extend it by linearity to the whole of $L(G)$. Here, $[g, h]$ is the usual commutator bracket in $G$, \ie \([g, h] = \inv{g} \inv{h} gh\).
\end{defin}
\begin{remark}
	For $i$-fold commutators, we work with left-normed commutators, \ie $[x_1, x_2, x_3] := [[x_1, x_2], x_3]$ and inductively, \([x_1, x_2, \ldots, x_n] := [[x_1, \ldots, x_{n - 1}], x_n].\)
\end{remark}
\begin{remark}
	We write cosets multiplicatively, \ie $g \gamma_{i}(G)$ and operations with cosets additively, \ie $g\gamma_{i}(G) + h \gamma_{i}(G)$.
\end{remark}

Each automorphism $\phi$ of $G$ induces an automorphism $\phi_*$ of $L(G)$ and if, moreover, each factor of the lower central series is finitely generated and torsion-free, we can talk about the eigenvalues of $\phi_*$: these are simply the eigenvalues of each $\phi_i$. The following result is then immediate from \cite[Lemma~2.2]{DekimpeGoncalves}

\begin{theorem}	\label{theo:usingL(G)ToEstablishRinf}
	Suppose $G$ is a finitely generated group such that all factors of the lower central series are torsion-free. Let $\phi \in \Aut(G)$ be an automorphism. Denote by $\phi_*$ the induced automorphism on $L(G)$ and by $\phi_i$ the restriction of $\phi_*$ to $L_i(G)$.
	
	If $\phi_i$ has eigenvalue $1$ for some $i$, then $R(\phi) = \infty$.

\end{theorem}

\section{Right-angled Artin groups}	\label{sec:RAAGs}
In this section, we briefly recall the necessary definitions and results regarding graphs and right-angled Artin groups, including isomorphisms between two related RAAGs and a generating set for the automorphism group of a RAAG.

\subsection{Definitions and examples}
 By a \emph{graph}, we mean a finite simple non-emtpy graph \(\Gamma(V, E)\) with vertices \(V\) and edges \(E\), although for technical reasons, we will sometimes mention the empty graph.

\begin{defin}
	Let $\Gamma(V, E)$ be a graph. The \emph{right-angled Artin group $A_\Gamma$ (or RAAG) associated to $\Gamma$} is defined by the presentation
	\[
		\pres{V}{[v, w], vw \in E}.
	\]
	The graph $\Gamma$ is called the \emph{defining graph of $A_{\Gamma}$}.
\end{defin}

From the definition of a RAAG, it is clear that the complete graph $K^n$ (\ie every two vertices are connected by an edge) on $n$ vertices corresponds to the free abelian group $\Z^n$ and that the edgeless graph $\overline{K^n}$ on $n$ vertices corresponds to the free group $F_n$.

\medskip
There are two operations on the level of graphs that correspond to natural operations on the level of groups and that will play an important role.
\begin{defin}
	Let $\Gamma_i(V_i, E_i)$, $i = 1, 2$, be two graphs. The \emph{disjoint union} of $\Gamma_1$ and $\Gamma_2$ is the graph
	\[
		\Gamma_1 \sqcup \Gamma_2(V_1 \sqcup V_2, E_1 \sqcup E_2)
	\]
	and the \emph{simplicial join} of $\Gamma_1$ and $\Gamma_2$ is the graph $\Gamma_1 * \Gamma_2$ with vertices $V_1 \sqcup V_2$ and edges
	\[
		E_{\Gamma_1 * \Gamma_2} = E_1 \sqcup E_2 \cup \{v_1v_2 \mid v_i \in V_i\}.
	\]
	We write $\sqcup_{n} \Gamma$ and $*_{n} \Gamma$ for the $n$-fold disjoint union, respectively, simplicial join of $\Gamma$ with itself.
\end{defin}
The following result follows readily from the definitions of a RAAG, direct product and free product.
\begin{prop}	\label{prop:directAndFreeProductRAAGs}
	Let $\Gamma_1, \Gamma_2$ be two graphs. Then $A_{\Gamma_1 \sqcup \Gamma_2} \cong A_{\Gamma_1} * A_{\Gamma_2}$ and $A_{\Gamma_1 * \Gamma_2} \cong A_{\Gamma_1} \times A_{\Gamma_2}$.
\end{prop}

\subsection{General isomorphisms of RAAGs}
In view of \cref{cor:characteristicQuotientRinf2}, it can be useful to transform one RAAG into another one by either deleting vertices or adding edges, and determining when this quotient is in fact characteristic. In this section, we make the first two claims more precise, in the next one, we discuss the characteristic quotients.

\begin{defin}
	Let $\Gamma(V, E)$ be a graph. We say that a subgraph $\Gamma'(V', E') \subseteq \Gamma$ is a \emph{full subgraph} or \emph{induced subgraph}, if $E'$ is given by
	\(
		\{vw \in E \mid v, w \in V'\}.
	\)
	
	Similarly, for $V' \subseteq V$, the \emph{subgraph induced on $V'$} is the graph $\Gamma'(V', E')$ where
	\(
		E' = \{vw \in E \mid v, w \in V'\}.
	\)
	We write $\Gamma(V')$.
\end{defin}

The RAAG associated to an induced subgraph $\Gamma(V')$ can be viewed as a subgroup of $A_\Gamma$ in a natural way.

\begin{lemma}[{\cite[Proposition~3.1]{DromsThesis}}]
	Let $\Gamma$ be a graph and $\Gamma' := \Gamma(V')$ an induced subgraph. The map
	\[
		i: A_{\Gamma'} \to A_{\Gamma}: v' \mapsto v', \quad \text{ for $v' \in V'$}
	\]
	is well-defined and injective.
\end{lemma}

We now formulate precisely how we can `delete vertices'. Given a subset \(S\) of a group \(G\), we denote the \emph{normal closure of \(S\) in \(G\)} by \(\normcl{S}_{G}\) or simply \(\normcl{S}\) if \(G\) is clear from the context.

\begin{prop}	\label{prop:eliminatingGenerators}
	Let $\Gamma$ be a graph, $\Gamma_1$ an induced subgraph and view $A_{\Gamma_1}$ as a subgroup of $A_\Gamma$. Let $\Gamma_2$ be the induced subgraph on $V_2 = V \setminus V_1$. We also write $\Gamma_2 = \Gamma \setminus \Gamma_1$. Define
	\[
		\phi: A_\Gamma \to A_{\Gamma_2}: v \mapsto \begin{cases}	1	&	\mbox{if } v \in V_1	\\
															v	&	\mbox{otherwise}
															\end{cases}.
	\]
	Then $\phi$ is a well-defined homomorphism, $\normcl{A_{\Gamma_1}}_{A_\Gamma} = \ker \phi$ and
	\[
		\frac{A_\Gamma}{\normcl{A_{\Gamma_1}}_{A_\Gamma}} \cong A_{\Gamma_2}.
	\]
\end{prop}
For a proof, we refer the reader to Fontelles master thesis \cite[Proposition~2.1]{FontellesThesis} or the master thesis \cite[Proposition~3.2.6]{SendenThesis} of one of the authors.

To make notations less heavy, we introduce the following definition:
\begin{defin}
	Let $\Gamma(V, E)$ be a graph and $A_\Gamma$ its associated RAAG. For a subset $V' \subseteq V$, the normal subgroup $\normcl{A_{\Gamma(V')}}$ is called the \emph{normal subgroup generated by $V'$} and we denote it by $N(V')$. Any subgroup $N(V')$ is called a \emph{normal vertex-subgroup}. 
\end{defin}
As we will never consider the subgroups $A_{\Gamma(V')}$ but only their normal closures, we simply refer to $N(V')$ as \emph{vertex-subgroups}.

Note that $N(V_1)N(V_2) = N(V_1 \cup V_2)$ for all $V_1, V_2 \subseteq V$. Indeed, both $N(V_1)$ and $N(V_2)$ lie in $N(V_1 \cup V_2)$, so their product does too. Conversely, it is clear that $A_{\Gamma(V_1 \cup V_2)}$ lies in $N(V_1)N(V_2)$, and as this product group is normal, we have that $N(V_1 \cup V_2) \leq N(V_1)N(V_2)$. 

The second isomorphism result for RAAGs tells us how we can `add edges'.
\begin{lemma}	\label{lem:imageOfNormalClosure}
	Let $G, H$ be groups, $S \subseteq G$ a subset and $\phi: G \to H$ a surjective homomorphism. Then $\phi(\normcl{S}) = \normcl{\phi(S)}$.
\end{lemma}

\begin{prop}	\label{prop:addingEdges}
	Let $\Gamma_1 \subseteq \Gamma_2$ be two graphs on $n$ vertices. Denote by $R_i$ the relators in the presentation of $A_{\Gamma_i}$ and put $N = \normcl{R_2}_{A_{\Gamma_1}}$. Then
	\[
		\frac{A_{\Gamma_1}}{N} \cong A_{\Gamma_2}.
	\]
\end{prop}
\begin{proof}
	Recall that $A_{\Gamma_1} \cong \frac{F_n}{\normcl{R_1}_{F_n}}$. Let $\pi: F_n \to A_{\Gamma_1}$ be the natural projection and put $\tilde{N} = \inv{\pi}(N)$. We claim that $\tilde{N} = \normcl{R_2}_{F_n}$. As $\normcl{R_1}_{F_n}$ lies in both subgroups, it is sufficient to prove that $\pi(\tilde{N}) = N$, by the correspondence theorem. This follows from \cref{lem:imageOfNormalClosure}, as $\pi$ is surjective. Hence, by the third isomorphism theorem, we have that
	\[
		\frac{A_{\Gamma_1}}{N} = \frac{F_n / \normcl{R_1}_{F_n}}{\normcl{R_2}_{F_n} / \normcl{R_1}_{F_n}} \cong \frac{F_n}{\normcl{R_2}_{F_n}} \cong A_{\Gamma_2}.	\qedhere
	\]
\end{proof}

\begin{cor}	\label{cor:abelianizationIsFreeAbelian}
	Let $\Gamma$ be a graph on $n$ vertices and $A_{\Gamma}$ its associated RAAG. Then $A_{\Gamma} / [A_{\Gamma}, A_{\Gamma}]$ is a free abelian group of rank $n$.
\end{cor}
\subsection{Generating set for \texorpdfstring{$\Aut(A_\Gamma)$}{}}
Naturally, to determine whether or not a RAAG has the $\Rinf$-property, we need to understand its automorphism group.

\begin{defin}
	Let $\Gamma(V, E)$ be a graph and $v \in V$. A vertex $w$ is called a \emph{neighbour of $v$}, or \emph{adjacent to $v$}, if $vw \in E$.
	
	We then define the \emph{link of $v$} as the set of all vertices adjacent to $v$, and it is denoted by $lk(v)$. The \emph{star of $v$} is $lk(v) \cup \{v\}$, and it is denoted by $st(v)$.	
	
	If $v \ne w$ are vertices, we say that \emph{$w$ dominates $v$} if $lk(v) \subseteq st(w)$, and write $v \leq w$.
\end{defin}
There are four basic types of automorphisms: \emph{graph automorphisms}, \emph{inversions}, \emph{transvections} and \emph{partial conjugations}.
\begin{itemize}
	\item Graph automorphisms of \(\Gamma\) extend to automorphisms of \(A_{\Gamma}\).
	\item Inversions \(\imath_{a}\) send one generator \(a\) to its inverse and leave all other generators fixed.
	\item Transvections are maps \(\tau_{ab}\) sending a generator \(a\) to \(ab\) and leaving all other generators fixed, where \(b\) is another (different) generator with \(a \leq b\).
	\item Partial conjugations are maps \(\gamma_{b, C}\) where \(b\) is a generator and \(C\) is a union of connected components of \(\Gamma \setminus \Gamma(st(b))\), sending every generator \(a\) in \(C\) to \(a^{b}\) and leaving the other generators fixed.
\end{itemize}

We will refer to these as \emph{automorphisms of basic type}. M.\ Laurence and H.\ Servatius have shown that the automorphisms of basic type generate $\Aut(A_\Gamma)$:
\begin{theorem}[\cite{Laurence95, Servatius89}]
	For every graph $\Gamma$, $\Aut(A_\Gamma)$ is generated by its graph automorphisms, inversions, transvections and partial conjugations.
\end{theorem}

\begin{lemma}		\label{lem:vertexSubgroupsArePartiallyCharacteristic}
	Let $\Gamma$ be a graph and $\Gamma_1$ an induced subgraph. Let ${\phi \in \Aut(A_\Gamma)}$ be a composition of inversions and partial conjugations. Then $\phi(\normcl{A_{\Gamma_1}}) \leq \normcl{A_{\Gamma_1}}$.
\end{lemma}
\begin{proof}
	Since \(\phi(N) \leq N\) and \(\psi(N) \leq N\) together imply \((\phi \circ \psi)(N) \leq N\) for any group \(G\) with normal subgroup \(N\) and automorphisms \(\phi\) and \(\psi\), it suffices to prove the statement for $\phi$ equal to an inversion or partial conjugation. If $\phi$ is an inversion, then $\phi(A_{\Gamma_1}) = A_{\Gamma_1}$, hence \(\phi(\normcl{A_{\Gamma_1}}) \leq \normcl{A_{\Gamma_1}}\). If $\phi$ is a partial conjugation, then $\phi(A_{\Gamma_{1}}) \subseteq \normcl{A_{\Gamma_{1}}}$ and thus \(\phi(\normcl{A_{\Gamma_1}}) \leq \normcl{A_{\Gamma_1}}\).
\end{proof}

So,  whenever we need to prove that $\normcl{A_{\Gamma_1}}$ is characteristic in $A_\Gamma$, we only need to prove it is preserved under graph automorphisms and transvections.
\section{Characteristic vertex-subgroups}	\label{sec:characteristicVertexSubgroups}

We begin the section by classifying all subsets of vertices that induce characteristic vertex-subgroups. After that, we treat two special cases, the vertices of maximal degree and the so-called transvection-free vertices, each of which implies a strong result regarding the $\Rinf$-property of RAAGs.
\subsection[Classification of characteristic vertex-subgroups]{Classification of characteristic vertex-subgroups}

For a graph $\Gamma(V, E)$ and a vertex $v \in V$, start with the set $V_0(v)$ containing only $v$. Add to $V_0(v)$ all vertices $w$ such that $\tau_{vw}$ is a well-defined transvection, to obtain $V_1(v)$. In symbols,
\[
	V_1(v) = \{w \in V \mid v \leq w\}.
\]
Inductively, given $V_i(v)$, put
\[
	V_{i + 1}(v) = \{w \in V \mid \exists v' \in V_i(v): v' \leq w\}.
\]
Define $V_\omega(v) = \bigcup_{i \in \N} V_i(v)$. Finally, put
\[
	\Vchar(v) = \bigcup_{\phi \in \Aut(\Gamma)} \phi(V_\omega(v)).
\]
\begin{defin}
	The set $\Vchar(v)$ defined above is called the \emph{characteristic closure of $v$ in $V$}.
\end{defin}
With this terminology, we then have the following result (which also explains the name \emph{characteristic closure}):
\begin{theorem}	\label{theo:classificationCharacteristicVertexSubgroups}
	Let $\Gamma$ be a graph and $A_\Gamma$ its associated RAAG. For vertices $v_1, \ldots, v_n \in V$, we have that
	\[
		N\left(\bigcup_{i = 1}^n \Vchar(v_i)\right)
	\]
	is characteristic in $A_\Gamma$.
	
	Conversely, if $N(V')$ is a characteristic vertex-subgroup for some $V' \subseteq V$, then $V'$ is a union of characteristic closures.
\end{theorem}
\begin{proof}
	As
	\[
		N\left(\bigcup_{i = 1}^n \Vchar(v_i)\right) = N(\Vchar(v_1)) \ldots N(\Vchar(v_n))
	\]
	and the product of characteristic subgroups is characteristic, it is sufficient to prove the theorem for $n = 1$. For simplicity, we write $V_{i}$ and $V_{\omega}$ for $V_{i}(v)$ and $V_{\omega}(v)$ as above.
	
	Let $\phi \in \Aut(A_\Gamma)$ be of basic type. By \cref{lem:vertexSubgroupsArePartiallyCharacteristic}, we only need to consider the case where $\phi$ is a graph automorphism or a transvection. If $\phi$ is a graph automorphism, then $\phi(w) \in \Vchar(v)$ for all $w \in \Vchar(v)$, by construction, hence $\phi(N(\Vchar(v))) \leq N(\Vchar(v))$. Similarly, if $\phi = \tau_{v'w}$ is a transvection with $v' \in \Vchar(v)$, then $v' = \psi(v'')$ for some $\psi \in \Aut(\Gamma)$ and $v'' \in V_\omega$. Note that $v'' \in V_i$ for some $i \in \N$. As $lk(v') \subseteq st(w)$, it follows that $lk(v'') \subseteq st(\inv{\psi}(w))$, hence $\inv{\psi}(w) \in V_{i + 1} \subseteq V_\omega$, and thus $w \in \psi(V_\omega) \subseteq \Vchar(v)$. It follows that $\phi(N(\Vchar(v))) \leq N(\Vchar(v))$. 
	
	Conversely, suppose $N(V')$ is characteristic in $A_\Gamma$ for some $V' \subseteq V$. For each $v \in V'$, we will prove that $\Vchar(v) \subseteq V'$; the result then follows from the equality
	\[
		V' = \bigcup_{v \in V'} \Vchar(v).
	\]
	First, we prove by induction on $i$ that for all $i \in \N$, $V_i$ as defined above is a subset of $V'$. Clearly, $V_0 \subseteq V'$. So suppose $V_i \subseteq V'$. Let $w \in V_{i + 1}$ be arbitrary. Then there is a $v' \in V_i$ such that $v' \leq w$. Then $\phi = \tau_{v'w}$ is a well-defined transvection, and as $N(V')$ is characteristic, $\phi(N(V')) \leq N(V')$. In particular, $v'w \in N(V')$ and as $v' \in V'$ by induction hypothesis, $w \in V'$. Consequently, $V_{i + 1} \subseteq V'$.
	
	From this, it follows that $V_\omega = \bigcup_{i \in \N} V_i \subseteq V'$. Now, let $\phi \in \Aut(\Gamma)$ and $v' \in V_\omega$. As $\phi(N(V')) \leq N(V')$, we see that $\phi(v') \subseteq V'$. Hence, $\phi(V_\omega) \subseteq V'$ for all $\phi \in \Aut(\Gamma)$. We conclude that
	\[
		\Vchar(v) = \bigcup_{\phi \in \Aut(\Gamma)} \phi(V_\omega) \subseteq V'.	\qedhere
	\]
\end{proof}

\subsection{Vertices of maximal degree}
In this section and the following, we will treat two particular instances of characteristic vertex-subgroups and each of them will have an important consequence in establishing the $\Rinf$-property of RAAGs. We begin by recalling the definition of the degree of a vertex.
\begin{defin}
	Let $\Gamma(V, E)$ be a graph and $v \in V$ a vertex. The \emph{degree of $v$}, denoted by $\deg(v)$, is the number of adjacent vertices of $v$.
	The \emph{maximal degree of $\Gamma$} is
	\[
		\Delta(\Gamma) := \max \{ \deg(v) \mid v \in V\}.
	\]
\end{defin}
\begin{lemma}
	Let $\Gamma$ be a graph and $v \in V$. Then for all $w \in \Vchar(v)$, $\deg(w) \geq \deg(v)$.
\end{lemma}
\begin{proof}
	First, we prove by induction that for all $i \in \N$ and $w \in V_i(v)$, we have $\deg(w) \geq \deg(v)$. Again, put $V_i = V_i(v)$. It clearly holds for $i = 0$. Suppose it holds for $i$. Let $w \in V_{i + 1}$. Then there is a $v' \in V_i$ such that $lk(v') \subseteq st(w)$. From this inclusion, it readily follows that $\deg(v') \leq \deg(w)$. As $\deg(v') \geq \deg(v)$ by induction hypothesis, we find that $\deg(w) \geq \deg(v)$.
	
	By taking the union over all $V_i$, it follows that $\deg(w) \geq \deg(v)$ for all $w \in V_\omega$. To show that it holds for all $w \in \Vchar(v)$, recall that
	\[
		\Vchar(v) = \bigcup_{\phi \in \Aut(\Gamma)} \phi(V_\omega)
	\]
	and note that every automorphism preserves the degree of a vertex, so if $w = \phi(w')$ for some $\phi \in \Aut(\Gamma)$ and $w' \in V_\omega$, then $\deg(w) = \deg(w') \geq \deg(v)$.
\end{proof}

\begin{cor}
	Let $\Gamma$ be a graph. Then $N(\Vmax)$ is characteristic in $A_\Gamma$, where
	\[
		\Vmax = \{v \in V \mid \deg(v) = \Delta(\Gamma)\}.
	\]
	Moreover, $A_\Gamma / N(\Vmax) \cong A_{\Gamma \setminus \Gamma(\Vmax)}$.
\end{cor}
\begin{proof}
	The first claim follows from \cref{theo:classificationCharacteristicVertexSubgroups} together with the equality
	\[
		\Vmax = \bigcup_{v \in \Vmax} \Vchar(v),
	\]
	which holds by the previous lemma. The isomorphism is a direct application of \cref{prop:eliminatingGenerators}.
\end{proof}

The true power of the previous corollary lies in the fact that for `almost all' graphs, $\Vmax$ will be non-trivial. Indeed, as each graph considered is assumed to be finite, $\Vmax \ne \emptyset$, and $\Vmax = V$ if and only if $\Gamma$ is \emph{regular}, \ie all vertices have the same degree. Therefore, every RAAG associated to a non-regular graph has a non-trivial characteristic vertex-subgroup, which provides a powerful tool in answering the conjecture. We elaborate further on this.

The next lemma is straightforward.

\begin{lemma}	\label{cor:characteristicQuotientsTransitive}
	Let $G$ be a group. If $G_{1}$ is a characteristic quotient of $G$ and $G_{2}$ one of $G_{1}$, then $G_{2}$ is a characteristic quotient of $G$.
\end{lemma}

\begin{lemma}[Simplification Lemma]	\label{lem:simplificationLemma}
	Let $\Gamma$ be a non-complete graph. Then $A_{\Gamma}$ has a RAAG $A_{\Gamma'}$ as a characteristic quotient where $\Gamma'$ is either
	\begin{itemize}
		\item disconnected;
		\item connected, regular and non-complete;
		\item non-regular, connected and such that $\Gamma'(V \setminus \Vmax)$ is complete.
	\end{itemize}
\end{lemma}
\begin{proof}
	Put $\Gamma_{0} = \Gamma$ and define inductively $\Gamma_{i + 1} = \Gamma_{i}(V_{i} \setminus V_{i, \mathrm{max}})$. Then each $A_{\Gamma_{i + 1}}$ is a characteristic quotient of $A_{\Gamma_{i}}$, and by \cref{cor:characteristicQuotientsTransitive}, also a characteristic quotient of $A_{\Gamma}$.
	
	Since $\Gamma$ is finite and for each $i$ we have $|V_{i + 1}| < |V_{i}|$, there is an $m$ such that $\Gamma_{m + 1}$ is the empty graph. This implies that $\Gamma_{m}$ is regular. If $\Gamma_{m}$ is disconnected, or connected and non-complete, we are done, so suppose $\Gamma_{m}$ is complete. Then $\Gamma_{m - 1}$ cannot be regular and $\Gamma_{m - 1}(V_{m - 1} \setminus V_{m - 1, \mathrm{max}})$ is complete. If $\Gamma_{m - 1}$ is disconnected, we are in the first case, otherwise we are in the third case.
\end{proof}

We give a name to the third type of graph arising in the Simplification Lemma.
\begin{defin}
	Let $\Gamma(V, E)$ be a connected non-regular graph. If the induced subgraph $\Gamma(V \setminus \Vmax)$ is complete, we call $\Gamma$ a \emph{(maximal degree)-by-(free abelian) graph}. We also say that $\Gamma$ is \emph{max-by-abelian}.
\end{defin}
\begin{remark}The inspiration for this terminology comes from group theory, where a group $G$ is called $\mathcal{P}$-by-$\mathcal{Q}$ if $G$ has a normal subgroup $N$ having property $\mathcal{P}$ such that $G / N$ has property $\mathcal{Q}$. Here, if $\Gamma$ is max-by-abelian, then the normal subgroup $N(\Vmax)$ of $A_\Gamma$ is generated by the vertices of \emph{maximal degree}, and the quotient $A_\Gamma / N(\Vmax)$ is a \emph{free abelian} group.

\end{remark}

The Simplification Lemma thus states that we only need to consider three types of graphs to establish the $\Rinf$-property for all non-abelian RAAGs. By demanding max-by-abelian graphs to be non-regular, there is no overlap between these three types of graphs. Note also that these types of graphs all have more structure than arbitrary graphs: although disconnectedness is not a very strong property from a graph theoretical point of view, it will be enough to ensure the $\Rinf$-property for the associated RAAG. For the second type, we have regularity, whereas the structure of max-by-abelian graphs is less obvious and less simple than that of the other two, as we will see later.

The Simplification Lemma proves the existence of at least one characteristic quotient of a particular form. It is of course possible for a RAAG to have multiple characteristic quotients, where for instance one quotient has a disconnected graph as defining graph and the other one a max-by-abelian graph.

\subsection{Transvection-free vertices}

Whereas $\Vmax$ provided us with the Simplification Lemma, the set of all so-called transvection-free vertices will be in some sense more constructive towards establishing the $\Rinf$-property for RAAGs: if the aforementioned set contains all vertices, then $A_\Gamma$ has the $\Rinf$-property.
\subsubsection{Definition and statement of the main theorem}
Let $\Gamma$ be a graph and $v \in V$ a vertex. In the construction of $\Vchar(v)$, we started with looking at which vertices $w \in V$ induce a well-defined transvection $\tau_{vw}$ on $A_\Gamma$. We consider now the vertices for which no such $w$ exist.
\begin{defin}
	Let $\Gamma(V, E)$ be a graph and $v \in V$ a vertex. We call $v$ \emph{transvection-free} if the set \(\{w \in V \mid v \leq w\}	\) only contains $v$. The set of all transvection-free vertices is denoted by $\Vnottrans$.
\end{defin}
\begin{remark}
	As $lk(v) \subseteq st(v)$, this definition makes sense.
\end{remark}

The following result is quite straightforward.
\begin{prop}	 \label{prop:transvectionFreeVerticesCharacteristic}
	For every graph $\Gamma$, the subgroup $N(\Vnottrans)$ is characteristic in $A_\Gamma$.
\end{prop}
\begin{proof}
	Let $v \in V$ be transvection-free. It immediately follows that $V_\omega(v) = \{v\}$, thus
	\[
		\Vchar(v) = \{\phi(v) \mid \phi \in \Aut(\Gamma)\}.
	\]
	If \(\phi \in \Aut(\Gamma)\) and $lk(\phi(v)) \subseteq st(w)$, then $lk(v) \subseteq st(\inv{\phi}(w))$, hence if $v$ is transvection-free, so is $\phi(v)$. Consequently,
	\[
		\Vnottrans = \bigcup_{v \in \Vnottrans} \Vchar(v)
	\]
	and we conclude by \cref{theo:classificationCharacteristicVertexSubgroups}.
\end{proof}

Although the previous proposition provides us with a characteristic subgroup of $A_\Gamma$ differing (in general) from $N(\Vmax)$, the most interesting situation occurs when $\Vnottrans = V$. In that case, no generator of $A_\Gamma$ admits transvections, so the group $\Aut(A_\Gamma)$ will be significantly smaller. In fact, it turns out that in that case, $A_\Gamma$ has the $\Rinf$-property, but even more is true:
\begin{restatable}{theorem}{transvectionFreeAutomorphism}	\label{theo:transvectionFreeAutomorphism}
	Let $\Gamma$ be a graph and denote by $\Autnottrans(A_{\Gamma})$ the subgroup of $\Aut(A_\Gamma)$ generated by all graph automorphisms, inversions and partial conjugations. If $\Gamma$ is not complete, then all $\phi \in \Autnottrans(A_{\Gamma})$ have infinite Reidemeister number, \ie $R(\phi) = \infty$ for all $\phi \in \Autnottrans(A_{\Gamma})$.
\end{restatable}

To prove this theorem, we need to understand the quotient groups $\gamma_i(A_\Gamma) / \gamma_{i + 1}(A_\Gamma)$ for $i = 1, 2, 3$. In order to do so, we study the associated Lie ring of a RAAG.

\subsubsection{Lyndon elements}
Lyndon elements can be used to show that the factors of the lower central series of a RAAG are free abelian groups, and they even provide us with a basis of these factors. This section is based on \cite{Wade}. Throughout this section, $\Gamma(V, E)$ is a graph, $V = \{v_1, \ldots, v_n\}$ and $A_\Gamma$ is its associated RAAG.

Define $W(V)$ to be the set of all words in the letters $v_1, \ldots, v_n$. The trivial word is denoted by $1$. The length of a word $w \in W(V)$ is denoted by $\len(w)$. For $w, w' \in W(V)$, we write $w \leftrightarrow w'$ if there are words $w_1, w_2 \in W(V)$ and vertices $v, v' \in V$ such that $vv' \in E$ and
\[
	w = w_1 vv' w_2	\quad \text{ and } \quad w' = w_1 v'v w_2.
\]
We then define an equivalence relation on $W(V)$: we write $w \sim w'$ if there are words $w_1, \ldots, w_k \in W(V)$ such that
\[
	w =  w_1 \leftrightarrow w_2 \leftrightarrow \ldots \leftrightarrow w_k = w'.
\]
It is clear that if $w \sim w'$ and $\tilde{w} \sim \tilde{w}'$, then $w\tilde{w} \sim w' \tilde{w}'$. This allows us to define a multiplication on the set $M := W(V) / \sim$. Also, the length of a word is preserved under the relation `$\leftrightarrow$', hence under the equivalence relation `$\sim$'. Thus, we can define the length of an element $m \in M$, which we still denote by $\len(m)$.

We put a total order on $W(V)$ as follows: for all $w \in W(V) \setminus \{1\}$, $1 < w$. Then, for $w \ne w' \in W(V) \setminus \{1\}$, write $w = v_i w_1$ and $w' = v_j w_2$ for $v_i, v_j \in V$ and $w_1, w_2 \in W(V)$. We put $w < w'$ if either $i < j$, or $i = j$ and $w_1 < w_2$.

We now proceed to define \emph{Lyndon elements} in $M$. Denote by $\pi: W(V) \to M$ the projection map linked to the equivalence relation defining $M$.
\begin{defin}
	For $m \in M$, we define the \emph{standard representative of $m$ in $W(V)$} to be the largest element of $\inv{\pi}(m)$ with respect to the total order. It is denoted by $std(m)$.
\end{defin}
Note that this is indeed well-defined: every element in $\inv{\pi}(m)$ has the same length, and as the alphabet $V$ is finite, so is the set of words of a fixed length. It follows that $\inv{\pi}(m)$ is finite, and as the order on $W(V)$ is total, this set has a largest element.

The standard representative of an element of $M$ allows us to define an order on $M$: for $m \ne m' \in M$, we put $m < m'$ if and only if $std(m) < std(m')$. From the totality of the order on $W(V)$, it follows that this order on $M$ is also total.

\begin{defin}
	Let $m, m' \in M$. We say that $m$ and $m'$ are \emph{transposed} if there exist $x, y \in M$ such that $m = xy$ and $m' = yx$.
\end{defin}
In $M$, being transposed is not an equivalence relation: if we take $\Gamma$ to be the graph with $V = \{v_1, v_2, v_3\}$ and $E = \{v_2v_3\}$, then putting $m_1 = v_2v_1v_3$, $m_2 = v_1v_3v_2 = v_1v_2v_3$ and $m_3 = v_3v_1v_2$, we find that $m_1$ and $m_2$ are transposed, as are $m_2$ and $m_3$, but $m_1$ cannot be transformed into $m_3$ using only one transposition.

To circumvent this, we do a similar trick as in the definition of $M$.
\begin{defin}
	Let $m, m' \in M$. We say that $m$ and $m'$ are \emph{conjugate} if there exist $m_1, \ldots, m_k$ such that $m = m_1$, $m' = m_k$, and $m_i$ and $m_{i + 1}$ are transposed for $1 \leq i \leq k - 1$.
\end{defin}
It is clear that this does define an equivalence relation on $M$.

\begin{defin}
	An element $m \in M$ is called \emph{primitive} if for every $x, y \in M$, the equality $m = xy = yx$ implies that $x = 1$ or $y = 1$.
	
	An element $m \in M$ is called a \emph{Lyndon element} if it is non-trivial, primitive and minimal in its conjugacy class \wrt the order $<$. 
\end{defin}
\begin{defin}
	Let $m \in M$. Then $init(m)$ is the set of all vertices in $V$ that can occur as the initial letter of a word in $\inv{\pi}(m)$.
\end{defin}
\begin{lemma}[{\cite[Corollary 3.2]{Lalonde}}]
	If $m \in M$ is a Lyndon element, then $init(m)$ consists of a single vertex.
\end{lemma}
For a Lyndon element $m$, we consider $init(m)$ to be the single vertex it contains rather than the singleton.
\begin{defin}
	For $m \in M$, we define $\zeta(m)$ to be
		\[
			supp(m) \cup \{v_j \mid \exists v_i \in supp(m): [v_i, v_j] \ne 1 \text{ in $A_\Gamma$}\}.
		\]
\end{defin}
\begin{remark}
	In the rest of this section, when writing $[v_i, v_j] \ne 1$, we mean $[v_i, v_j] \ne 1$ in $A_\Gamma$.
\end{remark}
Checking if a given element is a Lyndon element by means of the definition is quite cumbersome. Luckily, D.\ Krob and P.\ Lalonde \cite{Lalonde} proved the following equivalence.
\begin{theorem}[{\cite[Propositions 3.5 and 3.7]{Lalonde}}]	\label{theo:equivalenceLyndonElements}
	For $m \in M$, the following are equivalent:
	\begin{enumerate}[(i)]
		\item $m$ is a Lyndon element.
		\item For all $x, y \in M \setminus \{1\}$ such that $m = xy$, we have $m < y$.
		\item Either $\len(m) = 1$ or there exists Lyndon elements $x, y$ with $x < y$ and $init(y) \in \zeta(x)$ such that $m = xy$.
	\end{enumerate}
\end{theorem}

We finish by determining the Lyndon elements of length at most $3$.
\begin{prop}	\label{prop:LyndonElementsLength<=3}
	The sets
	\begin{align*}
		LE_1 &= \{v_i \mid 1 \leq i \leq n\}			\\
		LE_2 &= \{v_i v_j \mid 1 \leq i < j \leq n, [v_i, v_j] \ne 1\}		\\
		LE_3 &= \{v_i v_i v_k \mid 1 \leq i < k \leq n, [v_i, v_k] \ne 1\}	\\
			 &\quad \cup \{v_i v_j v_k \mid 1 \leq i < j < k \leq n, [v_i, v_j] \ne 1 \ne [v_i, v_k] \}	\\
			&  \quad \cup \{v_i v_j v_j \mid 1 \leq i < j \leq n, [v_i, v_j] \ne 1\}	\\
			& \quad \cup \{v_i v_j v_k \mid 1 \leq i < j \ne k \leq n, i < k, [v_i, v_j] \ne 1, ([v_i, v_k] \ne 1 \text{ or } [v_j, v_k] \ne 1)\}.
	\end{align*}
	are precisely the Lyndon elements of length $1$, $2$ and $3$, respectively.
\end{prop}
\begin{proof}
	For $LE_1$, this is clear. For $LE_2$, suppose $m$ is a Lyndon element of length $2$. Then $m = xy$ with $x < y$ both Lyndon elements and $init(y) \in \zeta(x)$. Note that $\len(x) = \len(y) = 1$, hence $x = v_i$ and $y = v_j$ for some $i, j$. As $x < y$, we have that $i < j$. Moreover, as $\zeta(x) = \{v_i\} \cup \{v_k \in V \mid [v_i, v_k] \ne 1\}$, we find that $[v_i, v_j] \ne 1$. Hence, $m = v_iv_j \in LE_2$. A similar argument shows that every element in $LE_2$ is a Lyndon element. 
	
	Suppose now that $m$ is a Lyndon element of length $3$. Write $m = v_iv_jv_k$. We distinguish two cases.
	
	\emph{Case 1: $v_i < v_jv_k$ are both Lyndon elements}. Then $v_jv_k \in LE_2$, hence $[v_j, v_k] \ne 1$ and $j < k$. As $v_i < v_jv_k$ and $std(v_i) = v_i$, $std(v_jv_k) = v_jv_k$, we find that $i \leq j$. Also, $v_j = init(v_jv_k) \in \zeta(v_i) = \{v_i\} \cup \{v_l \mid [v_i, v_l] \ne 1\}$. Hence, either $i = j$, or $i < j$ and $[v_i, v_j] \ne 1$.
	
	\emph{Case 2: $v_iv_j < v_k$ are both Lyndon elements}. Then $v_iv_j \in LE_2$, hence $i < j$ and $[v_i, v_j] \ne 1$. Similarly as in the previous case, $v_iv_j < v_k$ implies that $i < k$. Now, we have that $init(v_k) = v_k$ and
	\[
		\zeta(v_iv_j) = \{v_i, v_j\} \cup \{v_l \mid [v_i, v_l] \ne 1 \text{ or } [v_j, v_l] \ne 1\}.
	\]
	It follows that either
	\begin{itemize}
		\item $j = k$, or
		\item $j \ne k$ and $[v_i, v_k] \ne 1$, or
		\item $j \ne k$ and $[v_j, v_k] \ne 1$.
	\end{itemize}
	
	Conversely, suppose that $m = v_iv_jv_k \in LE_3$. In order to show that $m$ is a Lyndon element, we use the third criterion of \cref{theo:equivalenceLyndonElements} and case distinction.
	
	\emph{Case 1: $i = j < k$, $[v_i, v_k] \ne 1$:} Then both $v_i$ and $v_i v_k$ are Lyndon elements, $v_i < v_i v_k$ and $init(v_iv_k) = v_i \in \zeta(v_i)$, hence $m = v_iv_iv_k$ is a Lyndon element as well.
	
	\emph{Case 2: $i < j < k$, $[v_i, v_j] \ne 1 \ne [v_i, v_k]$:} In this case, $v_i$ and $v_jv_k$ are Lyndon elements, $v_i < v_jv_k$ and $init(v_jv_k) = v_j \in \zeta(v_i)$, as $[v_i, v_j] \ne 1$. Hence, $m = v_iv_jv_k$ is a Lyndon element as well.
	
	\emph{Case 3: $i < j = k$, $[v_i, v_j] \ne 1$:} As $[v_i, v_j] \ne 1$ and $i < j$, $v_iv_j$ is a Lyndon element, as is $v_j$. Since $std(v_iv_j) = v_iv_j$, it follows that $v_iv_j < v_j$ in $M$. Lastly, note that $init(v_j) = v_j \in \zeta(v_iv_j)$, so $m = v_iv_jv_j$ is a Lyndon element.
	
	\emph{Case 4: $i < j \ne k > i, [v_i, v_j] \ne 1$ and $([v_i, v_k] \ne 1$ or $[v_j, v_k] \ne 1)$:} Put $x = v_iv_j$ and $y = v_k$. As $i < j$ and $[v_i, v_j] \ne 1$, $x$ is a Lyndon element, and clearly, so is $y$. Note that $init(y) = v_k$ and $\zeta(x) = \{v_i, v_j\} \cup \{v_l \mid [v_i, v_l] \ne 1 \text{ or } [v_j, v_l] \ne 1\}$. As either $[v_i, v_k] \ne 1$ or $[v_j, v_k] \ne 1$, it follows that $init(y) \in \zeta(x)$ and therefore, $m$ is a Lyndon element.
\end{proof}

This explicit description of $LE_2$ and $LE_3$ will be particularly useful to determine linearly independent subsets of $\gamma_2(A_\Gamma) / \gamma_3(A_\Gamma)$ and $\gamma_3(A_\Gamma) / \gamma_4(A_\Gamma)$.

\medskip
Let $m$ be a Lyndon element of length at least $2$. By \cref{theo:equivalenceLyndonElements}, there exist Lyndon elements $x < y$ such that $m = xy$ and $init(y) \in \zeta(x)$. However, there can be multiple factorizations of this form. The factorization of $m$ where $y$ is minimal is called the \emph{standard factorization of $m$} and we write $S(m) = (x, y)$.

Using this standard factorization, we can define a bracketing procedure on the Lyndon elements: let $m$ be a Lyndon element. If $\len(m) = 1$, then $[m] = m$. If $\len(m) \geq 2$, write $(x, y) = S(m)$. Then the bracketing of $m$ is defined as $[[x], [y]]$, where $[x], [y]$ is the bracketing of $x$ and $y$, respectively. This bracketing can be interpreted in $A_{\Gamma}$ as the commutator bracket.

R.\ Wade \cite{Wade} then essentially proved the following:
\begin{theorem}	\label{theo:RAAGsAreFGTorsionFreeLCS}
	Let $\Gamma$ be a graph and $A_\Gamma$ its associated RAAG. Then the factors of the lower central series of $A_\Gamma$ are finitely generated and torsion-free. Moreover, after bracketing, the Lyndon elements of length $i$ form a $\Z$-basis of $L_i(A_\Gamma)$ for all $i$.
\end{theorem}

First of all, this theorem proves that we may use \cref{theo:usingL(G)ToEstablishRinf} for RAAGs. Secondly, combining this theorem with \cref{prop:LyndonElementsLength<=3}, we find explicit bases of $L_{2}(A_{\Gamma})$ and $L_{3}(A_{\Gamma})$. For $L_{2}(A_{\Gamma})$, this is immediate:
\begin{prop}	\label{prop:basisL2(AGamma)}
	The set
	\[
		\{[v_{i}, v_{j}]\gamma_{3}(A_{\Gamma}) \mid 1 \leq i < j \leq n, [v_{i}, v_{j}] \ne 1\}
	\]
	is a basis of $L_{2}(A_{\Gamma})$.
\end{prop}
Turning the Lyndon elements of length $3$ into an explicit basis of $L_{3}(A_{\Gamma})$ gives a quite long and ugly set. However, in the proof of \cref{theo:transvectionFreeAutomorphism}, we only need a linearly independent subset of this basis.
\begin{prop}	\label{prop:linearlyIndependentSubsetL3(AGamma)}
	Let $1 \leq i < j \leq n$ be indices such that $[v_{i}, v_{j}] \ne 1$. Then the bracketing of the Lyndon elements $v_{i}v_{i}v_{j}$ and $v_{i}v_{j}v_{j}$ is given by $[v_{i}, [v_{i}, v_{j}]]$ and $[[v_{i}, v_{j}], v_{j}]$, respectively.
	
	In particular, the set
	\begin{align*}
		&\{ [v_{i}, v_{j}, v_{i}]\gamma_{4}(A_{\Gamma}) \mid 1 \leq i < j \leq n, [v_{i}, v_{j}] \ne 1\}	\, \cup \\
		&\{ [v_{i}, v_{j}, v_{j}]\gamma_{4}(A_{\Gamma}) \mid 1 \leq i < j \leq n, [v_{i}, v_{j}] \ne 1\}
	\end{align*}
	forms a linearly independent subset of $L_{3}(A_{\Gamma})$.
\end{prop}
\begin{proof}
	\cref{prop:LyndonElementsLength<=3} implies that $v_{i}v_{i}v_{j}$ and $v_{i}v_{j}v_{j}$ are indeed Lyndon elements. Now, as neither $v_{i}v_{i}$ nor $v_{j}v_{j}$ is a Lyndon element, the standard factorizations are given by
	\[
		S(v_{i}v_{i}v_{j}) = (v_{i}, v_{i}v_{j}), \quad S(v_{i}v_{j}v_{j}) = (v_{i}v_{j}, v_{j}),
	\]
	hence the bracketing procedure gives $[v_{i}, [v_{i}, v_{j}]]$ and $[[v_{i}, v_{j}], v_{j}]$.
	
	The claim regarding the linearly independent subset follows from \cref{theo:RAAGsAreFGTorsionFreeLCS} and the fact that
	\[
		[v_{i}, [v_{i}, v_{j}]]\gamma_{4}(A_{\Gamma}) = -[v_{i}, v_{j}, v_{i}]\gamma_{4}(A_{\Gamma}). 	\qedhere
	\]
\end{proof}

\subsubsection{Proof of \cref{theo:transvectionFreeAutomorphism}}
We now have all the background needed to give a proof of \cref{theo:transvectionFreeAutomorphism}, which we restate here.

\transvectionFreeAutomorphism*

The main idea of the proof is to use \cref{theo:usingL(G)ToEstablishRinf}, \ie find for a given automorphism $\phi \in \Autnottrans(A_{\Gamma})$, an index $i$ such that the induced automorphisms $\phi_i$ on $L_i(A_\Gamma)$ has eigenvalue $1$. More precisely, the proof consists of the following steps:
\begin{enumerate}[Step 1.]
	\item We argue that it is sufficient to only consider the automorphisms in $\Autnottrans(A_{\Gamma})$ generated by graph automorphisms and inversions, by showing that partial conjugations descend to the trivial automorphism of $L(A_\Gamma)$.
	\item Denote by $\Autip(A_{\Gamma})$ the subgroup of $\Autnottrans(A_{\Gamma})$ generated by all graph automorphisms and inversions. We show that every element $\psi$ in $\Autip(A_{\Gamma})$ is conjugate to an element $\psi'$ of a certain `nice' form. By \cref{prop:conjugateEndomorphisms}, $R(\psi) = R(\psi')$.
	\item Finally, we show that each such `nice' $\psi'$ has infinite Reidemeister number, by finding an eigenvalue $1$ in either $L_1(A_\Gamma)$, $L_2(A_\Gamma)$ or $L_3(A_\Gamma)$.
\end{enumerate}
Throughout this section, we fix a non-complete graph $\Gamma$ on $n$ vertices and its associated RAAG $A_\Gamma$. We start with the first of two technical lemmata, which can be proven by straightforward induction.
\begin{lemma}	\label{lem:congruenceModuloLCS}
	Let $G$ be a group.
	\begin{enumerate}[(i)]
		\item If $x, y, z \in G$ are such that $x \equiv y \bmod \gamma_i(G)$, then $[x, z] \equiv [y, z] \bmod \gamma_{i + 1}(G)$.	\label{item:equalModGammai}
		\item \label{item:conjugationModGammai}For every $k \geq 1$, $x_1, \ldots, x_k \in G$ and $w_1, \ldots, w_k \in G$, we have that
		\[
			[x_1^{w_1}, \ldots, x_k^{w_k}] \equiv [x_1, \ldots, x_k] \bmod \gamma_{k + 1}(G).
		\]
		\item \label{item:exponentiationModGammai}For every $k \geq 1$, $x_1, \ldots, x_k \in G$ and $n_1, \ldots, n_k \in \Z$, we have that
		\[
			[x_1^{n_1}, \ldots, x_k^{n_k}] \equiv [x_1, \ldots, x_k]^{n_1 \ldots n_k} \bmod \gamma_{k + 1}(G).	
		\]
	\end{enumerate}
\end{lemma}

\begin{cor}	\label{cor:conjugationIsTrivialInLiering}
	For every partial conjugation $\phi \in \Aut(A_\Gamma)$ and every $i \geq 1$, we have that the induced automorphism $\phi_i$ on $\gamma_i(A_\Gamma) / \gamma_{i + 1}(A_\Gamma)$ is the identity map.
	
	Moreover, the induced automorphism $\phi_*$ on $L(A_\Gamma)$ is also the identity map.
\end{cor}
\begin{proof}
	As $\gamma_i(A_\Gamma) / \gamma_{i + 1}(A_\Gamma)$ is generated by the cosets of the $i$-fold commutators of the generators $a_1, \ldots, a_n$ of $A_\Gamma$, it is sufficient to show that $\phi_i$ is the identity map on these commutators. But since $\phi$ is given by conjugation for each generator, \cref{lem:congruenceModuloLCS} implies that
	\[
		\phi([a_{i_1}, \ldots, a_{i_i}]) = [\phi(a_{i_1}), \ldots, \phi(a_{i_i})] \equiv [a_{i_1}, \ldots, a_{i_i}] \bmod \gamma_{i + 1}(A_\Gamma).
	\]
	Hence, $\phi_i$ is the identity map. The `moreover' statement is immediate.
\end{proof}
The following lemma gives a more concrete description of a given element in $\Autip(A_{\Gamma})$.
\begin{lemma}		\label{lem:descriptionAutomorphismsTildeA}
	For $\phi \in \Autip(A_{\Gamma})$, there exist a permutation $\sigma \in S_n$ and integers $e_1, \ldots, e_n \in \{-1, 1\}$ such that
	\[
		\phi(a_i) = a_{\sigma(i)}^{e_i}
	\]
	for all $i \in \{1, \ldots, n\}$ and such that the map $a_i \mapsto a_{\sigma(i)}$ is a graph automorphism of $\Gamma$.
\end{lemma}
\begin{proof}
	The result trivially holds for the identity map, graph automorphisms and inversions, and it easily seen that the given form is preserved under composition. Hence, the result follows.
\end{proof}

The following lemma concerns the determinant of matrices of a specific form which will occur frequently in the proof of \cref{theo:transvectionFreeAutomorphism}.
\begin{lemma}		\label{lem:specialDeterminant}
	Let $k \geq 1$ be an integer and suppose $e_1, \ldots, e_k$ are all either equal to $1$ or $-1$. Define
	\[
		P(e_1, \ldots, e_k) := \begin{pmatrix}
			0		&	0		&	 \ldots	&	0		&	e_k	\\
			e_1		&	0		&	 \ldots	&	0		&	0	\\
			0		&	e_2		&	 \ldots	&	0		&	0	\\
			\vdots	&	\vdots	&	\ddots	&	\vdots	&	\vdots	\\
			0		&	0		&	 \ldots	&	0		&	0	\\
			0		&	0		&	 \ldots	&	e_{k - 1}	&	0.
		\end{pmatrix}
	\]
	Note that $P(e_1) = (e_1)$. Then
	\[
		\det(I_k - P(e_1, \ldots, e_k)) = 1 - \prod_{i = 1}^k e_i.
	\]
	In particular, $P(e_1, \ldots, e_k)$ has eigenvalue $1$ if and only if an even number of the $e_i$ equal $-1$.
\end{lemma}
\begin{proof}
	For $k = 1$, it is clear that $\det(I_1 - P(e_1)) = 1 - e_1$. For $k \geq 2$, expanding the determinant $\det(I_k - P(e_1, \ldots, e_k))$ along the first row yields the desired result.
\end{proof}

The last lemma (and the second technical one) makes more concrete what we meant by automorphisms in $\Autip(A_{\Gamma})$ of a `nice' form.
\begin{lemma}		\label{lem:conjugationInAtilde}
	Let $\phi \in \Autip(A_{\Gamma})$. Then there is a graph automorphism $\psi$ and an automorphism $\iota$, which is the composition of (possibly zero) inversions, such that $\phi$ and $\psi \circ \iota$ are conjugate.
	
	Moreover, if $\sigma \in S_n$ is such that $\psi(a_i) = a_{\sigma(i)}$ and $\sigma = c_1 \circ \ldots \circ c_k$ is the disjoint cycle decomposition of $\sigma$, we can choose $\iota$ such that each cycle $c_j = (c_{j, 1} \ldots c_{j, n_{j}})$ contains at most one number $i_j$ with $\iota(a_{i_j}) = \inv{a}_{i_j}$. (Here, we also write cycles of length $1$, e.g. the identity permutation is written as $(1)(2) \ldots (n)$).
\end{lemma}

\begin{proof}
	By \cref{lem:descriptionAutomorphismsTildeA}, we know that there is a permutation $\sigma \in S_n$ and elements $e_1, \ldots, e_n \in \{-1, 1\}$ such that $\phi(a_i) = a_{\sigma(i)}^{e_i}$ and such that $\psi: \Gamma \to \Gamma: a_i \mapsto a_{\sigma(i)}$ is a well-defined graph automorphism of $\Gamma$. Therefore, $\psi$ induces a graph automorphism of $A_\Gamma$, which we still denote by $\psi$.
	
	Write $\sigma = c_1 \circ \ldots \circ c_k$ in disjoint cycle decomposition. After renumbering the generators, we can assume that there are integers $1 = i_1 < i_2 < \ldots < i_k$ such that $c_j = (i_j \ldots i_{j + 1} - 1)$, where we put $i_{k + 1} - 1 = n$. If we consider $\phi_1: A_\Gamma / \gamma_2(A_\Gamma) \to A_\Gamma / \gamma_2(A_\Gamma)$, then the matrix of $\phi_1$ with respect to the basis $a_1, \ldots, a_n$ is of the form
	\[
		M := \begin{pmatrix}
			P(e_{i_1}, \ldots, e_{i_2 - 1})	&	0						&	0						&	\ldots	&	0						\\
			0						&	P(e_{i_2}, \ldots, e_{i_3 - 1})	&	0						&	\ldots	&	0						\\
			0						&	0						&	P(e_{i_3}, \ldots, e_{i_4 - 1})	&	\ldots	&	0						\\
			\vdots					&	\vdots					&	\vdots					&	\ddots	&	\vdots					\\
			0						&	0						&	0	 					&	\ldots	&	P(e_{i_k}, \ldots, e_n)
		\end{pmatrix}
	\]
	with each $P(e_{i_j}, \ldots, e_{i_{j + 1} - 1})$ as in \cref{lem:specialDeterminant}.

	If we can find a diagonal matrix $D$ with only $\pm 1$ on the diagonal such that $DMD$ is of the form
	\[
		\begin{pmatrix}
			P(\pm1, 1, \ldots, 1)	&	0				&	0				&	\ldots	&	0				\\
			0				&	P(\pm1, 1, \ldots, 1)	&	0				&	\ldots	&	0				\\
			0				&	0				&	P(\pm1, 1, \ldots, 1)	&	\ldots	&	0				\\
			\vdots			&	\vdots			&	\vdots			&	\ddots	&	\vdots			\\
			0				&	0				&	0	 			&	\ldots	&	P(\pm1, 1, \ldots, 1)
		\end{pmatrix}
	\]
	with each block of the same dimension as in $M$, we are done. Indeed, the matrix $D$ will then be the matrix of $\jmath_1: A_\Gamma / \gamma_2(A_\Gamma) \to A_\Gamma / \gamma_2(A_\Gamma)$ \wrt the basis $a_1, \ldots, a_n$, with $\jmath \in \Aut(A_\Gamma)$ a composition of inversions and $DMD$ corresponds to the automorphism $\psi \circ \iota$, where $\iota$ maps $a_{i_j}$ to $a^{\pm 1}_{i_j}$ for all $1 \leq j \leq k$ according to the sign in the $j$-th block of $DMD$ and leaves all other generators fixed. Then $\psi \circ \iota = \jmath \circ \phi \circ \jmath$ and $\iota$ will satisfy the `moreover'-part of the lemma, since each $P$-block contains at most one $-1$.
	
	It is sufficient to find such a diagonal matrix $D$ for each $P$-block and then putting all these blocks into one matrix. For ease of notation, we put $m = i_2 - 1$ and consider the block $P(e_1, \ldots, e_m)$. If $m = 1$, there is nothing to prove, so assume $m \geq 2$. Denote by $D_i$ the diagonal matrix with $1$'s on the diagonal except for the $i$-th position, there we put $-1$. Note that any product of the $D_i$'s is a diagonal matrix with only $\pm 1$ on the diagonal. It is not hard to see that
	\[
		D_iP(e_1, \ldots, e_m)D_i = P(e_1, \ldots, e_{i - 2}, -e_{i - 1}, -e_i, e_{i + 1}, \ldots, e_m),
	\]
	where $e_{0} = e_{m}$. Note that the parity of the number of $-1$'s is left unchanged. By starting with the $-1$ with the highest index and moving to $e_1$, we can clear out all $-1$'s, except for one if there were an odd number of them to begin with. Wherever this last $-1$ is situated, we can move it to the first position by conjugating with suitable $D_i$'s. In the end, we end up with $P(\pm 1, 1, \ldots, 1)$ and we are done.
\end{proof}
At last, we can give the proof of \cref{theo:transvectionFreeAutomorphism}.
\begin{proof}[Proof of \cref{theo:transvectionFreeAutomorphism}]
	Let $\phi \in A$. Since we will work with the induced morphism $\phi_*$ on $L(A_\Gamma)$, we may assume by \cref{cor:conjugationIsTrivialInLiering} that $\phi \in \Autip(A_{\Gamma})$. By \cref{lem:conjugationInAtilde}, $\phi$ is conjugate to $\psi \circ \iota$ with $\psi$ a graph automorphism and $\iota$ a composition of (possible zero) inversions satisfying the `moreover'-part of the statement of \cref{lem:conjugationInAtilde}. As $\phi$ and $\psi \circ \iota$ are conjugate, $R(\phi) = R(\psi \circ \iota)$ by \cref{prop:conjugateEndomorphisms}. So, we can assume that $\phi = \psi \circ \iota$.
	
	Let $\sigma \in S_n$ be the permutation associated to $\psi$ and write $\sigma = c_1 \circ \ldots \circ c_k$ in disjoint cycle decomposition. Denote by $\phi_i$ the induced automorphism on $L_i(A_\Gamma)$. By \cref{theo:RAAGsAreFGTorsionFreeLCS}, we can apply \cref{theo:usingL(G)ToEstablishRinf} to find that $R(\phi) = \infty$ if $\phi_i$ has eigenvalue $1$ for some $i$. We will distinguish several cases, and as each case will end with the sentence `$\phi_i$ has eigenvalue $1$', we will not mention \cref{theo:usingL(G)ToEstablishRinf} each time.
	
	\medskip	
	\emph{Case 1:} Suppose there is a cycle, say, $c_1$, such that $\iota$ does not invert any of the generators $a_i$ with $i \in c_1$. After renumbering, we can assume that $c_1 = (1 \ldots m)$ for some $m \geq 1$. Then, on $L_1(A_\Gamma)$, we have that
	\[
		\phi_1(a_1\ldots a_m\gamma_2(A_\Gamma)) = a_2a_3\ldots a_ma_1\gamma_2(A_\Gamma),
	\]
	so $\phi_1$ has a fixed point, \ie $1$ is an eigenvalue of $\phi_1$, hence $R(\phi) = \infty$.
	
	From now on, assume that each cycle contains an index $i_j$ such that $\iota(a_{i_j}) = \inv{a}_{i_j}$.
	
	\medskip	
	\emph{Case 2:} $k = n$, \ie each cycle in the decomposition of $\sigma$ has length $1$. Then $\phi(a_i) = \inv{a}_i$ for all $i$. As $\Gamma$ is not complete, there are $a_i \ne a_j$ with $a_ia_j \notin E$. Hence, $[a_i, a_j]\gamma_3(A_\Gamma)$ is non-trivial and
	\[
		\phi_2([a_i, a_j]\gamma_3(A_\Gamma)) = [\inv{a}_i, \inv{a}_j]\gamma_3(A_\Gamma) = [a_i, a_j]\gamma_3(A_\Gamma)
	\]
	by \cref{lem:congruenceModuloLCS}\eqref{item:exponentiationModGammai}. Then $\phi_2$ has eigenvalue $1$, and thus $R(\phi) = \infty$.
	
	From now on, assume that $k < n$.
	
	\medskip	
	\emph{Case 3:} There is a cycle, say, $c_1$ containing indices $i < j$ with $a_ia_j \notin E$. Again, we can assume that $c_1 = (1 \ldots m)$. As $\sigma$ induces a graph automorphism, we have that $a_{i + l}a_{j + l}$ is not an edge of $E$ either for all $1 \leq l \leq m$ (here, we work with indices modulo $m$ where we use \(1\) up to \(m\) as representatives, rather than \(0\) up to \(m - 1\)). After renumbering the generators, we can assume that $a_1$ is mapped onto $\inv{a}_2$. Consider then the set
	\[
		B := \{[a_{i + l}, a_{j + l}]\gamma_3(A_\Gamma) \mid 0 \leq l \leq m - 1\}
	\]
	Note that $B$ does not necessarily form a linearly independent set: if $i \equiv j + l \bmod m$ and $j \equiv i + l \bmod m$ for some $1 \leq l \leq m - 1$, then $[a_i, a_j] = - [a_{i + l}, a_{j + l}]$. These two congruences can be simultaneously fulfilled if and only if $2(j - i) \equiv 0 \bmod m$. If this condition is not satisfied, then $B$ is indeed a linearly independent set: each element in $B$ can be rewritten (up to sign) such that the first index is strictly less than the second. No two elements will have the same indices, and by \cref{prop:basisL2(AGamma)}, $B$ is a subset of a basis of $L_2(A_{\Gamma})$. We then can proceed as follows: first remark that
	\[
		\phi_2([a_{i + l}, a_{j + l}]\gamma_3(A_\Gamma)) = \begin{cases}
			-[a_{i + l + 1}, a_{j + l + 1}]\gamma_3(A_\Gamma)	&	\parbox[t]{.25\textwidth}{if $i + l \equiv 1 \bmod m$ or $j + l \equiv 1 \bmod m$}	\\[1em]
			[a_{i + l + 1}, a_{j + l + 1}]\gamma_3(A_\Gamma)	&	\mbox{otherwise}.
		\end{cases}
	\]
	Note that $i + l$ and $j + l$ cannot both be congruent to $1$ modulo $m$ simultaneously, as otherwise $i \equiv j \bmod m$, which is impossible as $1 \leq i < j \leq m$ (and thus also $m \geq 2$). As all elements in $B$ are distinct, there will be precisely two elements $[a_{i + l}, a_{j + l}]\gamma_3(A_\Gamma)$ that are mapped to $-[a_{i + l + 1}, a_{j + l + 1}]\gamma_3(A_\Gamma)$. Moreover, $\phi_2(\Span_\Z(B)) = \Span_\Z(B) =: V$, hence we can consider the matrix of $\phi_2$ restricted to $V$ with respect to this basis and find a matrix $P(e_1, \ldots, e_m)$ where precisely two of the $e_i$ are equal to $-1$ and the rest equals $1$. Hence, \cref{lem:specialDeterminant} implies that $P(e_1, \ldots, e_m)$ has eigenvalue $1$, so $\phi_2$ restricted to $V$ does too. This implies that $\phi_2$ has eigenvalue $1$, and consequently, $R(\phi) = \infty$.
	\medskip
	
	Now, suppose that $2(j - i) \equiv 0 \bmod m$. If $m = 2$, then $B$ contains $[a_1, a_2]\gamma_3(A_\Gamma)$ and in that case, 
	\[
		\phi_2([a_1, a_2]\gamma_3(A_\Gamma)) = -[a_2, a_1]\gamma_3(A_\Gamma) = [a_1, a_2]\gamma_3(A_\Gamma).
	\]
	Consequently, $\phi_2$ has eigenvalue $1$ and therefore $R(\phi) = \infty$.
	
	Finally, if $m > 2$ and $2(j - i) \equiv 0 \bmod m$, note that $0 < j - i < m$, hence $m$ is even and $j - i = m / 2$. The set 
	\begin{align*}
		\tilde{B}	&:= \{[a_{1 + l}, a_{m / 2 + 1 + l}, a_{1 + l}]\gamma_4(A_\Gamma) \mid 0 \leq l \leq m / 2 - 1\}	\\
				& \cup \{[a_{1 + l}, a_{m / 2 + 1 + l}, a_{m / 2 + 1 + l}]\gamma_{4}(A_{\Gamma}) \mid 0 \leq l \leq m / 2 - 1\}
	\end{align*}
	will be linearly independent, as all elements in $\tilde{B}$ are distinct and $\tilde{B}$ forms a subset of a linearly independent set of $L_3(A_\Gamma)$, by \cref{prop:linearlyIndependentSubsetL3(AGamma)}. Also note that
	\begin{align*}
		\phi_3([a_{1 + l}, a_{m / 2 + 1 + l}, a_{1 + l}]\gamma_4(A_\Gamma))	&= [a_{2 + l}, a_{m / 2 + 2 + l}, a_{2 + l}]\gamma_4(A_\Gamma)	\\
		\phi_3([a_{1 + l}, a_{m / 2 + 1 + l}, a_{m / 2 + 1 + l}]\gamma_4(A_\Gamma)) &= [a_{2 + l}, a_{m / 2 + 2 + l}, a_{m / 2 + 2 + l}]\gamma_4(A_\Gamma)
	\end{align*}
	for all $0 \leq l \leq m/2 - 2$ and that
	\begin{align*}
		\phi_{3}([a_{m / 2}, a_{m}, a_{m / 2}]\gamma_{4}(A_{\Gamma}))		&= [a_{m / 2 + 1}, a_{1}, a_{m / 2 + 1}]\gamma_{4}(A_{\Gamma})	\\
																&= -[a_{1}, a_{m / 2 + 1}, a_{m / 2 + 1}]\gamma_{4}(A_{\Gamma}),	\\
		\phi_{3}([a_{m / 2}, a_{m}, a_{m}]\gamma_{4}(A_{\Gamma}))		&= [a_{m / 2 + 1}, a_{1}, a_{1}]\gamma_{4}(A_{\Gamma})	\\
																&= -[a_{1}, a_{m / 2 + 1}, a_{1}]\gamma_{4}(A_{\Gamma})
	\end{align*} Hence, putting $W = \Span_\Z(\tilde{B})$, we find that $\phi_3(W) = W$ and that the matrix of $\phi_3$ restricted to $W$ with respect to $\tilde{B}$ is of the form $P(e_{1}, \ldots, e_{2m})$ as in \cref{lem:specialDeterminant}, where precisely two of the $e_{i}$'s are $-1$. This matrix has eigenvalue $1$, hence so does $\phi_3$. We conclude that $R(\phi) = \infty$.

	\medskip
	
	From now on, assume that each cycle $c_j$ satisfies the following property: if $c_j = (i_1 \ldots i_j)$, the induced subgraph $\Gamma(\{a_{i_1}, \ldots a_{i_j}\})$ is complete. As $\Gamma$ itself is not complete, there are cycles, say $c_1 = (1 \ldots m)$ and $c_2 = (m + 1 \ldots m + l)$, such that there are $i \in \{1, \ldots, m\}$, $j \in \{1, \ldots, l\}$ with $a_ia_{m + j} \notin E$. For notational convenience, we put $b_\imath= a_{\imath + m}$ for $\imath \in \{1, \ldots, l\}$. Hence, $a_i b_j \notin E$. It follows that also $a_1 b_{j - i + 1} \notin E$, as $\sigma^{1 - i}$ maps the non-edge $a_i b_j$ to the non-edge $a_{i + 1 - i} b_{j + 1 - i}$. Again, after renumbering, we can assume that $j - i + 1 = 1$. However, then it does not necessarily hold that $b_1 \mapsto \inv{b}_2$. By use of conjugation, we can obtain this nonetheless: suppose $b_\alpha \mapsto \inv{b}_{\alpha + 1}$ with $\alpha \ne 1$. Putting
	\(
		\jmath = \iota_2 \circ \iota_3 \circ \ldots \circ \iota_\alpha,
	\)
	where $\iota_k$ is the inversion $b_k$ to $\inv{b}_k$, a similar argument as in \cref{lem:conjugationInAtilde} yields that $\jmath \circ \phi \circ \jmath$ maps $b_1$ to $\inv{b}_2$, $b_\beta$ to $b_{\beta + 1}$ for $\beta \ne 1$ and coincides for the rest with $\phi$. Thus, we can assume that $a_1b_1 \notin E$.
	
	Put $K = \lcm(l, m)$ and
	\[
		B = \{[a_1, b_1]\gamma_3(A_\Gamma), \ldots, [a_K, b_K]\gamma_3(A_\Gamma) \}.
	\]
	Again, we consider the indices modulo $m$ (for $a_i$) and $l$ (for $b_j$), respectively. Then $B$ is a linearly independent subset of $L_2(A_\Gamma)$ and $V := \Span_\Z(B)$ satisfies $\phi_2(V) = V$. We count the number of elements in $B$ that are mapped to minus the next generator in $B$. For $1 \leq i \leq K$, we have that
	\[
		\phi_2([a_i, b_i]\gamma_3(A_\Gamma)) = \begin{cases}
			-[a_{i + 1}, b_{i + 1}]\gamma_3(A_\Gamma)	&	\parbox[t]{.4\textwidth}{if $i \equiv 1 \bmod m$ or $i \equiv 1 \bmod l$, but not both}	\\[1em]
			[a_{i + 1}, b_{i + 1}]\gamma_3(A_\Gamma)	&	\mbox{otherwise}.
		\end{cases}
	\]
	Note that $i \equiv 1 \bmod m$ and $i \equiv 1 \bmod l$ if and only if $i \equiv 1 \bmod K$, hence if and only if $i = 1$. The number of indices in $\{1, \ldots, K\}$ that are congruent to $1$ modulo $m$ is $K / m$, and similarly there are $K / l$ indices congruent $1$ modulo $l$. As we have counted $1$ twice, there are $K (\inv{m} + \inv{l}) - 1$ indices $i$ that are congruent to $1$ modulo $m$ or $l$. It follows that there are $K' := K (\inv{m} + \inv{l}) - 2$ elements in $B$ such that $\phi_2([a_i, b_i]\gamma_3(A_\Gamma)) = -[a_{i + 1}, b_{i + 1}]\gamma_3(A_\Gamma)$. 
	
	\emph{Case 4:} Suppose $K'$ is even. Then the matrix of $\phi_2$ restricted to $V$ with respect to the basis $B$ will be of the form $P(e_1, \ldots, e_{K'})$ with an even number of $e_i$ equal to $-1$. Applying \cref{lem:specialDeterminant} gives that $P(e_1, \ldots, e_{K'})$ and hence $\phi_2$ has eigenvalue $1$, so $R(\phi) = \infty$. 
	
	\emph{Case 5:} $K'$ is odd. Then either $K / l$ or $K / m$ is even, say, $K / m$ (the case $K / l$ even is analogous). The set
	\[
		B' := \{[a_1, b_1, b_1]\gamma_4(A_\Gamma), \ldots, [a_K, b_K, b_K]\gamma_4(A_\Gamma) \}.
	\]
	is a linearly independent subset of $L_3(A_\Gamma)$, by \cref{prop:linearlyIndependentSubsetL3(AGamma)}. Moreover, by \cref{lem:congruenceModuloLCS}\eqref{item:exponentiationModGammai}
	\[
		\phi_3([a_i, b_i, b_i]\gamma_4(A_\Gamma)) = \begin{cases}
			-[a_{i + 1}, b_{i + 1}, b_{i + 1}]\gamma_4(A_\Gamma)	&	\mbox{ if $i \equiv 1 \bmod m$}	\\
			[a_{i + 1}, b_{i + 1}, b_{i + 1}]\gamma_4(A_\Gamma)	&	\mbox{ otherwise}.
		\end{cases}
	\]
	Hence, precisely $K / m$ elements of $B'$ are mapped to minus the next element and $V' := \Span_\Z(B')$ satisfies $\phi_3(V') = V'$. The matrix of $\phi_3$ restricted to $V'$ with respect to $B'$ will then again be of the form as in \cref{lem:specialDeterminant} with an even number of $e_i$'s equal to $-1$, hence $\phi_3$ has eigenvalue $1$ and we can conclude that $R(\phi) = \infty$.
\end{proof}
\begin{defin}
	Let $\Gamma(V, E)$ be a graph. We call $\Gamma$ \emph{transvection-free} if $\Gamma$ is not $K^{1}$ and $\Vnottrans = V$.
\end{defin}
\begin{remark}
	For $K^{1}$, the sole vertex $v$ is, strictly speaking, transvection-free. Since $A_{K^{1}} = \Z$ is abelian, we exclude $K^{1}$ from the transvection-free graphs.
\end{remark}
The following important corollary is immediate.
\begin{cor}	\label{cor:transvectionFreeGraphsRinf}
	Let $\Gamma$ be a transvection-free graph. Then $A_\Gamma \in \Rinf$.
\end{cor}
\begin{proof}
	If $A_\Gamma$ is non-abelian and does not admit any transvections, then $\Autnottrans(A_{\Gamma}) = \Aut(A_\Gamma)$ and \cref{theo:transvectionFreeAutomorphism} implies that every $\phi \in \Aut(A_\Gamma) = \Autnottrans(A_{\Gamma})$ has infinite Reidemeister number.
\end{proof}

Of course, a theorem whose proof deserves a separate section should be of significant value. There are indeed a certain amount of transvection-free graphs, which we will discuss in the following section, but first we would like to make the following remark which explains to some extent the length of the proof: \cref{theo:transvectionFreeAutomorphism} holds for \emph{every} non-abelian RAAG. The only assumption regarding $\Gamma$ we needed in the proof was the fact that there are two vertices which are not connected by an edge. As only the complete graph does not meet this requirement, we consider a huge amount of very different graphs in the proof, which explains the number of cases.

The condition that $\Gamma$ is non-complete is also crucial: for $\Gamma = K^{n}$, we have $A_{\Gamma} = \Z^{n}$ and the automorphism $-{\Id_{\Z^{n}}}$ is a composition of graph automorphisms and inversions, but an easy calculation gives $R(-{\Id_{\Z^{n}}}) = 2^{n}$.

\subsubsection{Examples of transvection-free graphs}
To illustrate the significance of \cref{theo:transvectionFreeAutomorphism} and \cref{cor:transvectionFreeGraphsRinf}, we give an example of a family of transvection-free graphs. Later, in the section on regular graphs, we will see that most of the strongly regular graphs are also transvection-free, providing yet another example.
\begin{example}[Cycle graphs]	\label{ex:cycleGraphs}
	Let $n \geq 3$ be an integer. The \emph{cycle graph $C_n$} is the graph with vertex set $\{v_0, \ldots, v_{n - 1}\}$ and edge set \(\{v_i v_{i + 1} \mid i \in \{0, \ldots, n - 1\}\},\)
	where the indices are viewed modulo $n$. Note that $C_n$ is connected and regular for each $n \geq 3$, and non-complete for $n \geq 4$.
		\begin{figure}[h]
		\centering
		\begin{tikzpicture}
			\GraphInit[vstyle = Simple]
			\tikzset{VertexStyle/.style = {shape = circle,fill = black,minimum size = 4pt,inner sep=0pt}}
			\SetUpEdge[lw = 0.4pt]
			\grEmptyCycle[RA = 1.5, prefix = a, rotation = 10]{9}
			\EdgeInGraphLoop{a}{9}
		\end{tikzpicture}
		\caption{Cycle graph $C_9$}
	\end{figure}
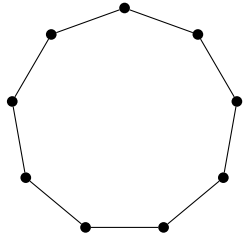
	
	We claim that $C_n$ is transvection-free if and only if $n \geq 5$. Suppose $n \geq 5$. For $i \in \{0, \ldots, n - 1\}$, we have that $lk(v_i) = \{v_{i - 1}, v_{i + 1}\}$ and $st(v_i) = \{v_{i - 1}, v_i, v_{i + 1}\}$. If $lk(v_i) \subseteq st(v_j)$, then either $v_i = v_j$, or $v_{i - 1} = v_{j + 1}$ and $v_{i + 1} = v_{j - 1}$. The latter case implies that $i - 1 \equiv j + 1 \bmod n$ and $i + 1 \equiv j - 1 \bmod n$. Hence, $-2 \equiv i - j \equiv 2 \bmod n$. As $n \geq 5$, this is impossible, hence $v_i = v_j$. We conclude that $C_n$ is indeed transvection-free. Consequently, $A_{C_n}$ has the $\Rinf$-property for $n \geq 5$.
	
	Conversely, if $n = 4$, note that $lk(v_1) = \{v_0, v_2\} = lk(v_3)$, so $C_4$ is not transvection-free. As $C_3 = K^3$, it is clear that also $C_3$ is not transvection-free.
\end{example}

\subsubsection{Properties of transvection-free graphs}
We end this section with some properties of transvection-free graphs.

\begin{defin}
	Let $\Gamma(V, E)$ be a graph. The \emph{complement graph $\cl{\Gamma}$} is the graph with vertex set $V$ and edge set
	\(
		\cl{E} = \{vw \mid v \ne w, vw \notin E\}.
	\)
	We also call $\cl{\Gamma}$ the \emph{complement of $\Gamma$}.
\end{defin}
For example, the complement of the complete graph $K^{n}$ is the edgeless graph $\cl{K^{n}}$, which also explains the notation. Moreover, it is not hard to see that
\(
	\cl{\Gamma_{1} * \Gamma_{2}} = \cl{\Gamma_{1}} \sqcup \cl{\Gamma_{2}}
\)
and that $\cl{\cl{\Gamma}} = \Gamma$.

The next result follows directly from the definitions.
\begin{lemma}	\label{lem:linksStarsComplement}
	Let $\Gamma$ be a graph and $v \in V$. Denote by $lk_{\Gamma}(v)$ the link of $v$ seen as vertex of $\Gamma$ and similarly for $st_{\Gamma}(v)$. Then $lk_{\cl{\Gamma}}(v) = V \setminus st_{\Gamma}(v)$ and $st_{\Gamma}(v) = V \setminus lk_{\cl{\Gamma}}(v)$.
\end{lemma}

\begin{defin}
	Let $\Gamma(V, E)$ be a graph. A vertex $v \in V$ is called \emph{isolated} if $lk(v) = \emptyset$.
\end{defin}
\begin{prop}	\label{prop:transvectionfreeGraphsUnderGraphsOperations}
	Let $\Gamma_{1}, \Gamma_{2}$ be transvection-free graphs. Then $\cl{\Gamma}_{1}$, $\Gamma_{1} \sqcup \Gamma_{2}$ and $\Gamma_{1} * \Gamma_{2}$ are all transvection-free as well.
\end{prop}
\begin{proof}
	A transvection-free graph cannot contain an isolated vertex $v$, since otherwise $lk(v) = \emptyset$ is contained in the star of every other vertex. Hence, the link of every vertex in $\Gamma_{1}$ and $\Gamma_{2}$ is non-empty, implying that any pair of dominating vertices in $\Gamma_{1} \sqcup \Gamma_{2}$ arises  in $\Gamma_{1}$ or $\Gamma_{2}$. Consequently, $\Gamma_{1} \sqcup \Gamma_{2}$ is transvection-free.
	\medskip
	
	For the complement of $\Gamma := \Gamma_{1}$: note that $w \in V$ dominates $v \in V$ seen as vertices in $\cl{\Gamma}$ if and only if $lk_{\cl{\Gamma}}(v) \subseteq st_{\cl{\Gamma}}(w)$. By \cref{lem:linksStarsComplement}, this is equivalent to $lk_{\Gamma}(w) \subseteq st_{\Gamma}(v)$, \ie $v$ dominates $w$ seen as vertices in $\Gamma$. As $\Gamma$ is transvection-free, it follows that $\cl{\Gamma}$ is transvection-free as well.
	\medskip
	
	Finally, consider the graph $\Gamma := \Gamma_{1} * \Gamma_{2}$. Suppose $v, w \in V_{1}$ are vertices with $v \leq w$ (in $\Gamma$). Since $V_{2} \subseteq lk(v), V_{2} \subseteq st(w)$, we have that $lk_{\Gamma}(v) = V_{2} \cup lk_{\Gamma_{1}}(v)$ and similarly for $st_{\Gamma}(w)$. This implies that $lk_{\Gamma_{1}}(v) \subseteq st_{\Gamma_{1}}(w)$, which contradicts the transvection-freeness of $\Gamma_{1}$. A similar argument holds if $v, w \in V_{2}$.
	
	If $v \in V_{1}$ and $w \in V_{2}$ are such that $lk_{\Gamma}(v) \subseteq st_{\Gamma}(w)$, then $V_{2} \subseteq lk_{\Gamma}(v) \subseteq st_{\Gamma}(w)$, hence, $st_{\Gamma_{2}}(w) = V_{2}$, which contradicts the fact that $\Gamma_{2}$ is transvection-free. We conclude that $\Gamma_{1} * \Gamma_{2}$ is transvection-free as well.
\end{proof}

\section{Disconnected graphs}	\label{sec:freeProductsRAAGs}
\addtocounter{subsection}{1}
From this section onwards, we discuss right-angled Artin groups associated to the three types of graphs arising from the Simplification Lemma. We start with disconnected graphs, which correspond to free products of RAAGs.

\begin{defin}
	A group \(G\) is called \emph{freely indecomposable} if \(G = G_{1} * G_{2}\) implies \(G_{1} = 1\) or \(G_{2} = 1\), \ie \(G\) cannot be written as a non-trivial free product.
\end{defin}
Clearly, abelian groups are freely indecomposable. Recently, D.\ Gon\c{c}alves, P.\ Sankaran and P.\ Wong proved the following \cite{GoncalvesSankaranWong}:
\begin{theorem}	\label{theo:RinfFreeProducts}
	Let \(n \geq 2\) and suppose \(G = G_{1}* \ldots* G_{n}\), where each \(G_{i}\) is freely indecomposable, and both \(G_{1}\) and \(G_{2}\) have a proper characteristic subgroup of finite index. Then \(G \in \Rinf\).
\end{theorem}

E.\ Green proved in her PhD-thesis \cite[Lemma 4.7]{Green} that a RAAG is freely indecomposable if and only if its defining graph is connected. Next, we show that each RAAG admits a proper characteristic subgroup of finite index.
\begin{lemma}	\label{lem:RAAGsAdmitCharacteristicSubgroup}
	Let \(\Gamma\) be a graph and \(A_{\Gamma}\) its associated RAAG. Then \(A_{\Gamma}\) has a proper characteristic subgroup of finite index.
\end{lemma}
\begin{proof}
	Suppose \(\Gamma\) has \(n\) vertices. Let \(p: A_{\Gamma} \to \frac{A_{\Gamma}}{[A_{\Gamma}, A_{\Gamma}]}\) be the natural projection. By \cref{cor:abelianizationIsFreeAbelian}, the abelianisation is isomorphic to \(\Z^{n}\). Note that \(2 \Z^{n}\) is a proper characteristic subgroup of finite index in \(\Z^{n}\). Then \(\inv{p}(2 \Z^{n})\) will be the desired subgroup of \(A_{\Gamma}\).
\end{proof}

We can therefore apply \cref{theo:RinfFreeProducts} to find the following result.
\begin{theorem}	\label{theo:RinfDisconnectedGraphs}
	Let \(\Gamma\) be a disconnected graph. Then \(A_{\Gamma}\) has the \(\Rinf\)-property.
\end{theorem}
\begin{proof}
	Write \(\Gamma = \bigsqcup_{i = 1}^{n} \Gamma_{i}\), where each \(\Gamma_{i}\) is connected and where \(n \geq 2\). By \cref{prop:directAndFreeProductRAAGs}, \(A_{\Gamma} \cong \Ast_{i = 1}^{n} A_{\Gamma_{i}}\). By \cite[Lemma 4.7]{Green}, each \(A_{\Gamma_{i}}\) is freely indecomposable and by \cref{lem:RAAGsAdmitCharacteristicSubgroup}, each \(A_{\Gamma_{i}}\) admits a proper finite index characteristic subgroup. Consequently, the conditions for \cref{theo:RinfFreeProducts} are fulfilled and we can conclude that \(A_{\Gamma} \in \Rinf\).
\end{proof}

\begin{remark}
	In  \cite[Chapter 5]{SendenThesis},  \cref{theo:RinfDisconnectedGraphs} was proved using the associated Lie ring of the RAAG.
\end{remark}
\section{Regular graphs}	\label{sec:regularGraphs}

We start this section with some results regarding direct products of RAAGs, which will then be applied to RAAGs associated to regular graphs. In particular, we discuss strongly regular graphs, which provide examples of both transvection-free graphs and direct products of RAAGs.

\subsection{Direct products of RAAGs}
For groups $G_{1}, \ldots, G_{n}$, the automorphism group of the direct product of the $G_{i}$'s always contains $\Aut(G_{1}) \times \ldots \times \Aut(G_{n})$, but the full automorphism group $\Aut(G_{1} \times \ldots \times G_{n})$ can be much bigger and much more difficult to describe. For RAAGs, however, a general description in terms of the factors is possible. We start with some general theory regarding automorphism groups of direct products before heading to direct products of RAAGs.

\subsubsection{Subgroups of \texorpdfstring{$\Aut(G_{1} \times \ldots \times G_{n})$}{}}	\label{subsec:subgroupsAutDirectProduct}
\begin{prop}	\label{prop:ReidemeisterNumberDirectProductAutomorphismGroups}
	Let $G_{1}, \ldots, G_{n}$ be groups and suppose (at least) one of them has the $\Rinf$-property. Then for each element $\phi$ of $\Aut(G_{1}) \times \ldots \times \Aut(G_{n}) \leq \Aut(G_{1} \times \ldots \times G_{n})$, we have $R(\phi) = \infty$.
\end{prop}
\begin{proof}
	Write $\phi = (\phi_{1}, \ldots, \phi_{n})$. It is clear that, for $(g_{1}, \ldots, g_{n}), (h_{1}, \ldots, h_{n}) \in G_{1} \times \ldots \times G_{n}$, we have
	\[
		(g_{1}, \ldots, g_{n}) \conj{\phi} (h_{1}, \ldots, h_{n}) \iff \forall i \in \{1, \ldots, n\}: g_{i} \conj{\phi_{i}} h_{i}.
	\]
	Thus,
	\[
		R(\phi) = \prod_{i = 1}^{n} R(\phi_{i}),
	\]
	from which the result immediately follows.
\end{proof}
For an \(n\)-fold direct product of a group with itself, we can consider a bigger subgroup of the automorphism group than merely the direct product: there is an injective group homomorphism
	\[
		\psi: \Aut(G) \wr S_{n} \to \Aut(G^{n}): (\phi, \sigma) = (\phi_{1}, \ldots, \phi_{n}, \sigma) \mapsto \psi(\phi, \sigma)
	\]
	where
	\[
		\psi(\phi, \sigma)(g_{1}, \ldots, g_{n}) = (\phi_{1}(g_{\inv{\sigma}(1)}), \ldots, \phi_{n}(g_{\inv{\sigma}(n)})).
	\]
Here, \(\Aut(G) \wr S_{n}\) is the wreath product, \ie the semidirect product \(\Aut(G)^{n} \rtimes S_{n}\), where the action is given by
\[
	\sigma \cdot (\phi_{1}, \ldots, \phi_{n}) = (\phi_{\inv{\sigma}(1)}, \ldots, \phi_{\inv{\sigma}(n)}).
\]

We will write elements in \(\Aut(G) \wr S_{n}\) simply as \((\phi, \sigma)\).

\begin{prop}	\label{prop:ReidemeisterNumberWreathProduct}
	Let $G$ be a group having the $\Rinf$-property and $n \geq 1$ an integer. Then for each $\chi \in \Aut(G) \wr S_{n} \leq \Aut(G^{n})$, we have $R(\chi) = \infty$.
\end{prop}
\begin{proof}
	Write $\chi = (\phi, \sigma)$. Suppose $(g, 1, \ldots, 1) \conj{\chi} (h, 1, \ldots, 1)$. Then there are $x_{1}, \ldots, x_{n} \in G$ with
	\begin{align}
		g &= x_{1}h\inv{\phi_{1}(x_{\inv{\sigma}(1)})}	\nonumber\\
		1 &= x_{i} \inv{\phi_{i}(x_{\inv{\sigma}(i)})} \quad \text{ $i \geq 2$}.	\label{eq:xi=phi(xi)}
	\end{align}
	First, suppose that $\inv{\sigma}(1) = 1$. Then $g \conj{\phi_{1}} h$, hence the map
	\[
		F: \{[(g, 1, \ldots, 1)]_{\chi} \mid g \in G\} \to \Reid[\phi_{1}]: [(g, 1, \ldots, 1)]_{\chi} \mapsto [g]_{\phi_{1}}
	\]
	is well-defined and surjective. The domain of $F$ is a subset of $\Reid[\chi]$, hence $R(\chi) \geq R(\phi_{1}) = \infty$, since $\phi_{1} \in \Aut(G)$ and $G \in \Rinf$.
	
	If $\inv{\sigma}(1) \ne 1$, let $m > 0$ be the smallest integer such that $\sigma^{-m}(1) = 1$. Repeatedly using \eqref{eq:xi=phi(xi)} gives
	\[
		g = x_{1}h\inv{\phi_{1}(x_{\inv{\sigma}(1)})} = x_{1}h\inv{\phi_{1}(\phi_{\inv{\sigma}(1)}(x_{\sigma^{-2}(1)}))} = \ldots = x_{1}h\inv{\tilde{\phi}(x_{1})}
	\]
	where
	\[
		\tilde{\phi} = \phi_{1} \circ \phi_{\inv{\sigma}(1)} \circ \ldots \circ \phi_{\sigma^{1 - m}(1)}.
	\]
	Note that $\tilde{\phi}$ only depends on $\chi$. Thus, $g \conj{\tilde{\phi}} h$ and the map
	\[
		F: \{[(g, 1, \ldots, 1)]_{\chi} \mid g \in G\} \to \Reid[\tilde{\phi}]: [(g, 1, \ldots, 1)]_{\chi} \mapsto [g]_{\tilde{\phi}}
	\]
	is well-defined and surjective. Similarly as before, we find that $R(\chi) = \infty$.
\end{proof}

\begin{cor}	\label{cor:wreathProductEqualAutomorphismGroupRinf}
	If a group \(G\) has the \(\Rinf\)-property and \(n \geq 1\) is an integer such that \(\Aut(G^{n}) = \Aut(G) \wr S_{n}\), then \(G^{n} \in \Rinf\).
\end{cor}

\subsubsection{Automorphism group of direct product of RAAGs}

The description of the automorphism group of a direct product of RAAGs is due to N.~Fullarton \cite{Fullarton}, and G.~Giovanni and N.~Wahl \cite{Gandini}, whose results we present in this section, together with its implications for the $\Rinf$-property for RAAGs. Recall that a direct product on group theoretical level corresponds to a simplicial join on graph theoretical level.

\begin{prop}[{\cite[Proposition 3.1]{Gandini}}]	\label{prop:primeDecompositionRAAGs}
	Let $\Gamma$ be a graph. Then $A_{\Gamma}$ admits a unique maximal decomposition as
	\[
		A_{\Gamma} = A_{\Gamma_{1}} \times \ldots \times A_{\Gamma_{k}}
	\]
	for induced subgraphs $\Gamma_{1}, \ldots, \Gamma_{k}$. This decomposition is unique up to isomorphism and permutation of the factors.
\end{prop}

\begin{lemma}[{\cite[Proposition 2.2]{CharneyVogtmann}}]	\label{lem:centreOfRAAG}
	Let $\Gamma$ be a graph and $A_{\Gamma}$ its associated RAAG. Then the centre of $A_{\Gamma}$ is given by
	\(
		Z(A_{\Gamma}) = \gen{\{v \in V \mid \deg(v) = |V| - 1\}}.
	\)
\end{lemma}

\begin{prop}[{\cite[Proposition 3.3]{Gandini}}]	\label{prop:furtherDecompositionAutomorphismGroupRAAG}
	Let $\Gamma$ be a graph with maximal decomposition
	\begin{equation}	\label{eq:simplicialDecompositionGamma}
		\Gamma = K^{d} * \Ast_{j = 1}^{k} (*_{i_{j}} \Gamma_{j})
	\end{equation}
	where all $\Gamma_{j}$ are non-isomorphic and non-complete graphs.
	Then
	\[
		\Aut(A_{\Gamma}) \cong \Z^{d |V'|} \rtimes (\GL_{d}(\Z) \times (\Aut(A_{\Gamma_{1}}) \wr S_{i_{1}}) \times \ldots \times (\Aut(A_{\Gamma_{k}}) \wr S_{i_{k}})).
	\]
\end{prop}
Using this description and the results from \cref{subsec:subgroupsAutDirectProduct}, we can describe what happens on the level of Reidemeister numbers.

\begin{theorem}	\label{theo:directProductsRAAGsRinf}
	Let $\Gamma$ be a graph with maximal decomposition
	\begin{equation*}
		\Gamma = K^{d} * \Ast_{j = 1}^{k} (*_{i_{j}} \Gamma_{j})
	\end{equation*}
	where all $\Gamma_{j}$ are non-isomorphic and non-complete graphs. If any graph $\Gamma_{j}$ is such that $A_{\Gamma_{j}} \in \Rinf$, then $A_{\Gamma} \in \Rinf$.
\end{theorem}
\begin{proof}
	Since $Z(A_{\Gamma}) = A_{K^{d}}$ by \cref{lem:centreOfRAAG} and the centre of a group is characteristic, we have by \cref{prop:eliminatingGenerators} the isomorphism
	\[
		\frac{A_{\Gamma}}{A_{K^{d}}} \cong A_{\Gamma'}
	\]
	where
	\begin{equation}	\label{eq:decompositionGamma'}
		\Gamma' = \Ast_{j = 1}^{k} (*_{i_{j}} \Gamma_{j}).
	\end{equation}
	The maximality of the decomposition of $\Gamma$ implies that \eqref{eq:decompositionGamma'} is the maximal decomposition of $\Gamma'$. Hence, \cref{prop:furtherDecompositionAutomorphismGroupRAAG} implies that
	\[
		\Aut(A_{\Gamma'}) \cong \Times_{j = 1}^{k} (\Aut(A_{\Gamma_{j}}) \wr S_{i_{j}}).
	\]
	Suppose that $A_{\Gamma_{j}} \in \Rinf$ for some $1 \leq j \leq k$. Then the isomorphism above combined with \cref{prop:ReidemeisterNumberDirectProductAutomorphismGroups,prop:ReidemeisterNumberWreathProduct} implies that $A_{\Gamma'} \in \Rinf$ as well. Since $A_{\Gamma'}$ is a characteristic quotient of $A_{\Gamma}$, we consequently have that $A_{\Gamma} \in \Rinf$.
\end{proof}

\begin{example}	\label{ex:finiteProductsStartingWithZ}
	As an application of the previous theorem, consider the smallest family $\G$ of groups satisfying the following two properties: $\G$ contains $\Z$ and $\G$ is closed under taking finite free products and finite direct products. We claim that any non-abelian group in $\G$ has the $\Rinf$-property. Note that $\G$ contains only RAAGs.

	Let $G \in \G$ be a non-abelian group and write $G = A_{\Gamma}$ for some (non-complete) graph $\Gamma$. The group $G$ splits as either a free or direct product of two groups in $\G$. In the former case, $\Gamma$ is disconnected, hence $G \in \Rinf$ by \cref{theo:RinfDisconnectedGraphs}. In the latter case, write $\Gamma = \Gamma_{1} * \ldots * \Gamma_{k}$ in maximal decomposition, with $k \geq 2$. Since $G$ is non-abelian, one factor, say, $\Gamma_{1}$, is non-complete. As the decomposition of $\Gamma$ is maximal and $A_{\Gamma_{1}} \in \G$, we must have that $\Gamma_{1}$ is disconnected. Consequently, $A_{\Gamma_{1}} \in \Rinf$ by \cref{theo:RinfDisconnectedGraphs}. We can then apply \cref{theo:directProductsRAAGsRinf} to conclude.
	
	In particular, any (non-abelian) finite direct product of free groups has the $\Rinf$-property.
\end{example}

We would like to remark that \cref{theo:directProductsRAAGsRinf} does \emph{not} state that any direct product of RAAGs has the $\Rinf$-property if (at least) one factor has it. It does state so for certain types of RAAGs, namely those for which their corresponding graph does not split as a simplicial join. Nonetheless, the result is strong enough for our purposes.

\subsection{Regular graphs}
We formally state the definition of a regular graph.
\begin{defin}
	Let $\Gamma(V, E)$ be a graph and $k \geq 0$. We call $\Gamma$ \emph{$k$-regular} if $\deg(v) = k$ for all $v \in V$. We call $\Gamma$ \emph{regular} if there is some $k \geq 0$ such that $\Gamma$ is $k$-regular.
\end{defin}
For example, the complete graph $K^n$ is $(n - 1)$-regular, whereas all cycle graphs on at least $3$ vertices are $2$-regular. Clearly, there is (up to isomorphism) only one $0$-regular graph on $n$ vertices, namely $\cl{K^{n}}$. We start with some very basic properties of regular graphs. For a proof, we refer the reader to \cite[Chapter 2]{ChartrandZhang}.
\begin{lemma}[Handshaking Lemma]	\label{lem:handshakingLemma}
	For a graph $\Gamma$ the following equality holds:
	\[
		2|E| = \sum_{v \in V} \deg(v).
	\]
\end{lemma}
\begin{lemma}	\label{lem:basicPropertiesRegularGraphs}
	Let $\Gamma$ be a $k$-regular graph on $n$ vertices. Then
	\begin{enumerate}[(i)]
		\item $kn \equiv 0 \bmod 2$.				\label{item:knequiv0regularGraph}
		\item $\cl{\Gamma}$ is $(n - k - 1)$-regular.	\label{item:complementRegularGraph}
	\end{enumerate}
\end{lemma}
With this lemma, we classify all $k$-regular graphs on $n$ vertices with $k \in \{1, 2, n - 2, n - 3\}$.
\begin{prop}	\label{prop:1regularGraphs}
	Suppose $\Gamma$ is a $1$-regular graph on $n$ vertices. Then $\Gamma$ is the disjoint union of $n / 2$ copies of $K^{2}$.
\end{prop}
\begin{proof}
	Note that $n$ is even, since $n\cdot 1 \equiv 0 \bmod 2$ by \cref{lem:basicPropertiesRegularGraphs}\eqref{item:knequiv0regularGraph}. Write $\Gamma = \Gamma_{1} \sqcup \ldots \sqcup \Gamma_{k}$ with each $\Gamma_{i}$ connected. Then each $\Gamma_{i}$ is a connected $1$-regular graph, so $\Gamma_{i}$ contains at least $2$ vertices, say $v$ and $w$, that are connected by an edge. However, $\Gamma_{i}$ cannot contain more than $2$ vertices, since no additional vertex can be adjacent to either $v$ or $w$. Hence, $\Gamma_{i}$ is (isomorphic to) $K^{2}$, implying that $\Gamma$ is a disjoint union of $n / 2$ copies of $K^{2}$.
\end{proof}
Taking complements, we obtain all $(n - 2)$-regular graphs on $n$ vertices.
\begin{cor}	\label{cor:n-2RegularGraphs}
	Suppose $\Gamma$ is an $(n - 2)$-regular graph on $n$ vertices. Then $\Gamma$ is the simplicial join of $n / 2$ copies of $\cl{K^{2}}$.
\end{cor}
\begin{cor}	\label{cor:1/n-2RegularGraphsRinf}
	Suppose $\Gamma$ is a non-complete graph on $n$ vertices that is either $1$-regular or $(n - 2)$-regular. Then $A_{\Gamma} \in \Rinf$.
\end{cor}
\begin{proof}
	If $\Gamma$ is $1$-regular and non-complete, it is disconnected by \cref{prop:1regularGraphs}, hence $A_{\Gamma} \in \Rinf$ by \cref{theo:RinfDisconnectedGraphs}. If $\Gamma$ is $(n - 2)$-regular, it is a simplicial join of $n / 2$ copies of $\cl{K^{2}}$. Then \cref{ex:finiteProductsStartingWithZ} implies that $A_{\Gamma} \in \Rinf$.
\end{proof}

\begin{prop}	\label{prop:2regularGraphs}
	Suppose $\Gamma$ is a $2$-regular graph. Then $\Gamma$ is the disjoint union of cycle graphs.
\end{prop}
\begin{proof}
	It is sufficient to prove that every \emph{connected} $2$-regular graph is a cycle graph. So, suppose that $\Gamma$ is connected. Denote by $v_{1}, \ldots, v_{n}$ the vertices of $\Gamma$. Consider a path of maximal length $k$ in $\Gamma$. Without loss of generality, this is the path $v_{1}, v_{2}, \ldots, v_{k}$. For each $j \in \{2, \ldots, k - 1\}$, the link  of $v_{j}$ is given by $lk(v_{j}) = \{v_{j - 1}, v_{j + 1}\}$. Since $\deg(v_{k}) = 2$, there is a vertex besides $v_{k - 1}$ adjacent to $v_{k}$, say $v_{i}$. Note that $i \in \{1, k + 1, \ldots, n\}$. If $i \geq k + 1$, then the path was not maximal. Hence, $i = 1$ and since $\Gamma$ is connected, $k$ must be equal to $n$; otherwise, there would not exist a path from $v_{k + 1}$ to $v_{k}$. We conclude that $\Gamma$ is indeed a cycle graph.
\end{proof}
Again, taking complements gives us all $(n - 3)$-regular graphs on $n$ vertices.
\begin{cor}	\label{cor:n-3RegularGraphs}
	Suppose $\Gamma$ is an $(n - 3)$-regular graph on $n$ vertices. Then $\Gamma$ is the simplicial join of complements of cycle graphs.
\end{cor}

\begin{cor}	\label{cor:2/n-3RegularGraphsRinf}
	Suppose $\Gamma$ is a non-complete graph on $n$ vertices that is either $2$-regular or $(n - 3)$-regular. Then $A_{\Gamma} \in \Rinf$.
\end{cor}
\begin{proof}
	If $\Gamma$ is $2$-regular, then \cref{prop:2regularGraphs} implies $\Gamma$ is either disconnected or a cycle graph. In the former case, $A_{\Gamma} \in \Rinf$. In the latter, $n \geq 4$ since $\Gamma$ is non-complete. If $n = 4$, then $\Gamma = \cl{K^{2}} * \cl{K^{2}}$. Consequently, \cref{ex:finiteProductsStartingWithZ} implies that $A_{\Gamma} \in \Rinf$. If $n \geq 5$, then $\Gamma$ is transvection-free by \cref{ex:cycleGraphs}, so \cref{cor:transvectionFreeGraphsRinf} implies that $A_{\Gamma} \in \Rinf$.
	\medskip

	Now suppose that $\Gamma$ is $(n - 3)$-regular. Then $\Gamma = \cl{C}_{n_{1}} * \ldots * \cl{C}_{n_{k}}$ for some $k \geq 1$ and $n_{i} \geq 3$. Note that the complement of the cycle graph $C_{3}$ is $\cl{K^{3}}$, the complement of $C_{4}$ is the disjoint union $K^{2} \sqcup K^{2}$, whereas the complement of the cycle graph $C_{n}$ for $n \geq 5$ is transvection-free by \cref{ex:cycleGraphs} and \cref{prop:transvectionfreeGraphsUnderGraphsOperations}. Consequently, $A_{\cl{C}_{m}} \in \Rinf$ for all $m \geq 3$ by \cref{theo:RinfDisconnectedGraphs,cor:transvectionFreeGraphsRinf}. Moreover, as each cycle graph is connected, their complements do not split as simplicial joins. Hence, the maximal decomposition of $\Gamma$ is given by $\cl{C}_{n_{1}} * \ldots * \cl{C}_{n_{k}}$. Applying \cref{theo:directProductsRAAGsRinf} then gives the result.
\end{proof}

\subsection{Strongly regular graphs}
We follow \cite[Chapter 10]{GodsilRoyle} for the definition of strongly regular graphs. As the name suggests, these graphs have strong regularity conditions.
\begin{defin}
	Let $\Gamma(V, E)$ be a graph. We say that $\Gamma$ is \emph{strongly regular with parameters $n$, $k$, $\lambda$ and $\mu$} if \(\Gamma\) is a \(k\)-regular graph on \(n\) vertices such that
	\begin{itemize}
		\item every two adjacent vertices have $\lambda$ neighbours in common, i.e. if $vw \in E$, then $|lk(v) \cap lk(w)| = \lambda$,
		\item every two non-adjacent vertices have $\mu$ neighbours in common, i.e. if $vw \notin E$ and $v \ne w$, then $|lk(v) \cap lk(w)| = \mu$,
		\item $1 \leq k < n - 1$.
	\end{itemize}
	We also say that $\Gamma$ is an \emph{$srg(n, k, \lambda, \mu)$}.
\end{defin}
\begin{remark}
	The only $0$-regular graph on $n$ vertices is the edgeless graph $\overline{K^n}$. If we would consider this graph to be strongly regular, then $\lambda$ would be undefined. On the other hand, the only $(n - 1)$-regular graph on $n$ vertices is the complete graph $K^n$, for which the parameter $\mu$ is undefined. This explains why we require that $1 \leq k < n - 1$. The corresponding RAAGs of the aforementioned graphs are the free group and the free abelian group, respectively, and for both we already determined the Reidemeister spectrum. Hence, it is no loss to exclude these graphs from the strongly regular ones.
\end{remark}
Note that we thus can assume that $n \geq 2$ and in that case, $\lambda \leq k - 1$ and $\mu \leq k$. For an example of a strongly regular graph, we first need a definition.

\begin{defin}
	Let $\Gamma(V, E)$ be a graph and $p \geq 2$. We say that $\Gamma$ is \emph{$p$-partite} if $V$ admits a partition $V_1, \ldots, V_p$ such that no vertices in $V_i$ are connected by an edge. If $\Gamma$ is $p$-partite for some $p$, we say that $\Gamma$ is \emph{multipartite}.
\end{defin}
\begin{example}
	Let $p \geq 2$ and $n_1, \ldots, n_p$ be strictly positive integers. The \emph{complete $p$-partite graph of order $n_1, \ldots, n_p$} is the graph
	\[
		K(n_1, \ldots, n_p) := \overline{K^{n_1}} * \overline{K^{n_2}} *\ldots * \overline{K^{n_p}}.
	\]
	If $n_1 = n_2 = \ldots = n_p =: n$, we also write $K^p_n$. If $n = 1$, then $K^p_1$ is the complete graph on $p$ vertices, hence the notation coincides.
	
	\medskip
	Now, we claim that for $n \geq 2$, $K_n^p$ is an $srg(np, n(p - 1), n(p - 2), n(p - 1))$. Clearly, the number of vertices is $np$. If we denote by $V_1, \ldots, V_p$ the partition of the vertex set $V$, then $|V_i| = n$ for all $i$ and each $v \in V_i$ is connected with every vertex in $V \setminus V_i$. Hence, $\deg(v) = (p - 1)n$, so $K_n^p$ is $(p - 1)n$-regular.
	
	Two vertices that are not adjacent lie in the same $V_i$, hence they have $n(p - 1)$ common neighbours. If two vertices $v, w$ are adjacent, then $v \in V_i$ and $w \in V_j$ with $i \ne j$. A vertex $v' \in V$ is a neighbour of $v$ if and only if $v' \notin V_i$ and it is a neighbour of $w$ if and only if $v' \notin V_j$. Hence, the common neighbours of $v$ and $w$ are precisely all vertices in $V \setminus (V_i \cup V_j)$, which contains $n(p - 2)$ elements.
\end{example}
The regularity conditions for a strongly regular graph $\Gamma$ restrict the possibilities for characteristic vertex-subgroups: as $\Gamma$ is regular, $\Vmax = V$, and $\Vnottrans$ will be either empty or the whole of $V$. In order to prove this last statement, it will be more convenient to work with the complement of $\Vnottrans$, which we now give a name.
\begin{defin}
	Let $\Gamma(V, E)$ be a graph. We call the complement of $\Vnottrans$ the set of all \emph{transvection-admitting vertices} and we denote it by $\Vtrans$. Note that
	\[
		\Vtrans = \{v \in V \mid \exists w \ne v \in V: v \leq w\}.
	\]
\end{defin}

\begin{lemma}	\label{lem:symmetryVtransRegular}
	Let $\Gamma$ be a $k$-regular graph and $v, w \in V$ two vertices. Then
	\begin{enumerate}[(i)]
		\item $lk(v) \subseteq lk(w) \iff lk(w) \subseteq lk(v) \iff |lk(v) \cap lk(w)| = k$.
		\item $st(v) \subseteq st(w) \iff st(w) \subseteq st(v) \iff |st(v) \cap st(w)| = k + 1$.
	\end{enumerate}
\end{lemma}
\begin{proof}
	Both statements follow directly from the fact that $|lk(v)| = |lk(w)| = k$ and $|st(v)| = |st(w)| = k + 1$.
\end{proof}

\begin{prop}	\label{prop:transvectionfreeVerticesSRG}
	Let $\Gamma$ be an $srg(n, k, \lambda, \mu)$. Then
	\begin{equation}	\label{eq:VtransInSRG}
		\Vtrans = \begin{cases}
			\emptyset 	&	\mbox{ if $\lambda < k - 1$ and $\mu < k$}	\\
			V			&	\mbox{ otherwise}.
		\end{cases}
	\end{equation}
\end{prop}
\begin{proof}
	As mentioned before, we can assume that $n \geq 2$. Let $v, w \in V$ be distinct vertices. We look for equivalent conditions for $v \leq w$. We distinguish two cases.
	
	\emph{Case 1: $vw \notin E$}. Then $v \leq w$ is equivalent with $lk(v) \subseteq lk(w)$ and hence with $|lk(v) \cap lk(w)| = k$, by the previous lemma. As $\Gamma$ is strongly regular, $|lk(v) \cap lk(w)| = \mu$. Hence, $lk(v) \subseteq lk(w)$ if and only if $k = \mu$.
	
	\emph{Case 2: $vw \in E$}. In this case, $v \leq w$ is equivalent with $st(v) \subseteq st(w)$ and hence with $|st(v) \cap st(w)| = k + 1$. Note that
	\[
		|st(v) \cap st(w)| = |\{v, w\} \cup (lk(v) \cap lk(w))| = 2 + \lambda
	\]
	Therefore, $st(v) \subseteq st(w)$ if and only if $\lambda = k - 1$.
	
	\medskip
	With this information, we can prove \eqref{eq:VtransInSRG}: if $\lambda < k - 1$ and $\mu < k$, then $v \leq w$ can never happen, hence $\Vtrans = \emptyset$.
	
	If $\mu = k$, recall that $\Gamma$ is a $k$-regular non-complete graph. Hence, for every vertex $v$ there is a non-adjacent vertex $w$. Then $lk(v) \subseteq lk(w)$ as $\mu = k$. Consequently, $\Vtrans = V$.
	
	If $\lambda = k - 1$, take a vertex $v$ and one of its neighbours $w$ (which exists, as $k \geq 1$). As $\lambda = k - 1$, $st(v) \subseteq st(w)$. This holds for all $v \in V$, so $\Vtrans = V$.
\end{proof}

\begin{cor}	\label{cor:almostAllSRGhaveRinf}
	If $\Gamma$ is an $srg(n, k, \lambda, \mu)$ with $\lambda < k - 1$ and $\mu < k$, then $A_\Gamma \in \Rinf$.
\end{cor}
\begin{proof}
	From the previous proposition, it follows that $\Vtrans = \emptyset$, hence $\Gamma$ is transvection-free. Applying \cref{cor:transvectionFreeGraphsRinf} now yields the result.
\end{proof}
A natural question is to ask which strongly regular graphs have either $\lambda = k - 1$ or $\mu = k$. To answer this question, we first point out a relation between those two kinds of graphs.

\begin{lemma}[{\cite[p.~218]{GodsilRoyle}}]	\label{lem:complementSRGisSRG}
	Let $\Gamma$ be an $srg(n, k, \lambda, \mu)$. Then $\cl{\Gamma}$ is an $srg(n, n - k - 1, n - 2k - 2 + \mu, n - 2k + \lambda)$.
\end{lemma}
\begin{lemma}[{\cite[p.~219]{GodsilRoyle}}]	\label{lem:relationParametersSRG}
	Let $\Gamma$ be an $srg(n, k, \lambda, \mu)$. Then
	\[
		(n - k - 1) \mu = k (k - \lambda - 1).
	\]
	In particular, $\mu = k$ if and only if $\lambda = 2k - n$, and $\lambda = k - 1$ if and only if $\mu = 0$.
\end{lemma}

Hence, if $\Gamma$ is an $srg(n, k, k - 1, 0)$, then $\cl{\Gamma}$ is an $srg(n, k', \lambda', \mu')$ with
\begin{align*}
	k' &= n - k - 1	\\
	\lambda' &= n - 2k - 2 = 2k' - n	\\
	\mu' &= n - k - 1 = k'.
\end{align*}
So, we only have to determine the strongly regular graphs with $\mu = k$ in order to determine all srg's with either $\mu = k$ or $\lambda = k - 1$.

\begin{prop}[{\cite[Lemma 10.1.1]{GodsilRoyle}}] \label{prop:SRGwithlambda = k - 1}
	Let \(\Gamma\) be an \(srg(n, k, \lambda, \mu)\). If \(\lambda = k - 1\), then \(\Gamma\) is the disjoint union of \(\frac{n}{k + 1}\) copies of \(K^{k + 1}\).
\end{prop}
\begin{cor}	\label{cor:SRGwithmu = k}
	Let $\Gamma$ be an $srg(n, k, \lambda, \mu)$. If $\mu = k$, then $\Gamma$ is the complete multipartite graph $K^{n / (n - k)}_{n - k}$.
\end{cor}
\begin{proof}
	By \cref{lem:complementSRGisSRG} and \cref{prop:SRGwithlambda = k - 1}, \(\cl{\Gamma}\) is the disjoint union of \(\frac{n}{n - k}\) copies of \(K^{n - k}\). Taking complements again implies that \(\Gamma\) is the simplicial join of \(\frac{n}{n - k}\) copies of \(\cl{K^{n - k}}\), \ie \(\Gamma\) is the complete multipartite graph \(K^{n / (n - k)}_{n - k}\).
\end{proof}

The RAAGs associated to srg's with $\lambda = k - 1$ and $\mu = k$ are
\[
	\Ast_{i = 1}^{\frac{n}{k + 1}} \Z^{k + 1} \quad \text{ and } \quad \Times_{i = 1}^{\frac{n}{n - k}} F_{n - k},
\]
respectively. The first group has the $\Rinf$-property by \cref{theo:RinfDisconnectedGraphs}. For the second, we invoke \cref{ex:finiteProductsStartingWithZ}. Combining with \cref{cor:almostAllSRGhaveRinf}, we find the following result.
\begin{theorem}	\label{theo:SRGsHaveRinf}
	Let $\Gamma$ be a strongly regular graph. Then $A_{\Gamma} \in \Rinf$.
\end{theorem}
\section{Max-by-abelian graphs}	\label{sec:MBA}
Max-by-abelian graphs possess more structure than arbitrary (connected) graphs, but this structure is much less apparent than the structure of regular graphs, for instance. We therefore start by discussing general graph theoretical properties of max-by-abelian graphs, and we also provide some examples. Thereafter, we move on to their associated RAAGs and use the aforementioned structural properties to prove the $\Rinf$-property in certain cases.
\subsection{Examples and general properties}
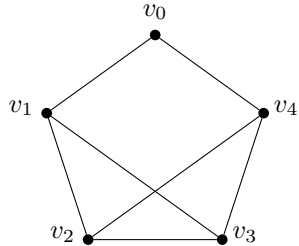
\begin{figure}[h]
	\centering
	\begin{tikzpicture}
		\GraphInit[vstyle = Simple]
		\SetVertexMath
		\tikzset{VertexStyle/.style = {shape = circle,fill = black,minimum size = 4pt,inner sep=0pt}}
	        \SetUpEdge[lw  = 0.4pt]  
		\grEmptyCycle[RA=1.5, prefix=v, rotation = 90]{5}%
		\tikzset{AssignStyle/.append style = {above= 2pt}}
		\AssignVertexLabel[color = black, size = \small]{v}{$v_{0}$, , , ,}
		\tikzset{AssignStyle/.append style = {left = 1pt}}
		\AssignVertexLabel[color = black, size = \small]{v}{, $v_{1}$, $v_{2}$, , }
		\tikzset{AssignStyle/.append style = {right = 1pt}}
		\AssignVertexLabel[color = black, size = \small]{v}{, , , $v_{3}$, $v_{4}$}
		\EdgeInGraphLoop{v}{5}
		\EdgeFromOneToSel{v}{v}{1}{3}
		\EdgeFromOneToSel{v}{v}{2}{4}
	\end{tikzpicture}
	\caption{A $(5, 4, 3)$-max-by-abelian graph}
	\label{fig:example543MBAgraph}
\end{figure}

We extend the definition of a max-by-abelian graph with some parameters.
\begin{defin}
	Let $\Gamma(V, E)$ be a graph. We say that $\Gamma$ is \emph{$(n, k, d)$-max-by-abelian} (or $(n, k, d)$-MBA) if
	\begin{itemize}
		\item $\Gamma$ is non-regular and connected,
		\item $\Gamma(V \setminus \Vmax)$ is complete,
		\item $|V| = n$, $|\Vmax| = k$ and $\Delta(\Gamma) = d$.
	\end{itemize}
\end{defin}
In particular, $k \leq n - 1$ and $d \leq n - 2$, since $\Gamma$ is non-regular (and thus a fortiori non-complete). \Cref{fig:example543MBAgraph} shows an example of a $(5, 4, 3)$-MBA graph with $\Vmax = \{v_{1}, v_{2}, v_{3}, v_{4}\}$.
\begin{example}	\label{ex:simplicialJoinDisjointUnionCompleteGraphs}
	Let $n, n_{1}, n_{2}$ be strictly positive integers with $n \leq n_{1}, n \leq n_{2}$ and $n^{2} < n_{1} n_{2}$. We claim that $\Gamma := (K^{n} \sqcup K^{n_{1}}) * (K^{n} \sqcup K^{n_{2}})$ is an MBA-graph. It is clearly connected. Denote the vertices of the first copy of $K^{n}$ by $a_{i}$, of the second copy by $c_{i}$, the vertices of $K^{n_{1}}$ by $b_{i}$ and those of $K^{n_{2}}$ by $d_{i}$. Then
	\begin{align*}
		\deg(a_{i})	&= n - 1 + n + n_{2} 	= 2n + n_{2} - 1	\\
		\deg(b_{i})	&= n_{1} - 1 + n + n_{2}	= n + n_{1} + n_{2} - 1	\\
		\deg(c_{i})	&= n - 1 + n + n_{1}		= 2n + n_{1} - 1	\\
		\deg(d_{i})	&= n_{2} - 1 + n + n_{1}	= n + n_{1} + n_{2} - 1
	\end{align*}
	Since at least one of the inequalities $n < n_{1}$ and $n < n_{2}$ holds, $\Gamma$ is not regular. We may assume that $n_{1} \leq n_{2}$. If $n < n_{1}$, then
	\[
		\Vmax = \{b_{i} \mid 1 \leq i \leq n_{1}\} \cup \{d_{i} \mid 1 \leq i \leq n_{2}\}
	\]
	and $\Gamma(V \setminus \Vmax) = K^{n} * K^{n} = K^{2n}$ is complete. In this case, $\Gamma$ is an $(2n + n_{1} + n_{2}, n_{1} + n_{2}, n + n_{1} + n_{2} - 1)$-MBA graph.
	
	If $n = n_{1}$, then $n < n_{2}$ and
	\(
		\Vmax = V \setminus \{c_{i} \mid 1 \leq i \leq n\}.
	\)
	In this case, $\Gamma(V \setminus \Vmax) = K^{n}$, so $\Gamma$ is an $(2n + n_{1} + n_{2}, n + n_{1} + n_{2}, n + n_{1} + n_{2} - 1)$-MBA graph.
	
	Note that the corresponding RAAG is given by $(\Z^{n} * \Z^{n_{1}}) \times (\Z^{n} * \Z^{n_{2}})$, which has the $\Rinf$-property by \cref{theo:RinfDisconnectedGraphs,theo:directProductsRAAGsRinf}.
\end{example}
The first properties we derive are inequalities providing relations amongst or restrictions on the parameters $n, k$ and $d$.
\begin{prop}
	Let $\Gamma$ be an $(n, k, d)$-MBA graph. Then
	\begin{equation}	\label{eq:boundsEdgesMBA}
		\binom{n}{2} - k(n - d - 1) \leq |E| \leq \frac{n(d - 1) + k}{2}.
	\end{equation}
\end{prop}
\begin{proof}
	The inequalities follow from a count of the minimal resp.\ maximal number of edges in $\Gamma$. First, we prove the upper bound. Every vertex in $\Vmax$ has degree $d$, so every vertex in $V \setminus \Vmax$ has at most degree $d - 1$. We thus have that
	\begin{equation*}
		|E| = \frac{1}{2} \sum_{v \in V} \deg(v) \leq \frac{k d + (n - k)(d - 1)}{2} = \frac{n(d - 1) + k}{2}.
	\end{equation*}
	The lower bound is more involved. Since the induced subgraph $\Gamma(V \setminus \Vmax)$ is the complete graph on $n - k$ vertices, $\Gamma$ has at least $\binom{n - k}{2}$ edges. We now count how many edges we deleted going from $\Gamma$ to $\Gamma(V \setminus \Vmax)$. Put $\Vmax = \{v_{1}, \ldots, v_{k}\}$. Deleting $v_{1}$ results in deleting $d$ edges, one for every neighbour of $v_{1}$. Deleting $v_{i}$ after deleting $v_{1}$ up to $v_{i - 1}$ results in deleting at least $d - i + 1$ edges. Indeed, the only edges with endpoint $v_{i}$ that could already have been deleted are those with $v_{j}$ with $1 \leq j \leq i - 1$ as (other) endpoint. There are at most $i - 1$ such vertices and since $\deg(v_{i}) = d$, there are at least $d - i + 1$ edges with one endpoint equal to $v_{i}$ and one in $V \setminus \Vmax$. Hence, deleting $v_{i}$ after deleting $v_{1}$ up to $v_{i - 1}$ results in deleting at least $d - i + 1$ edges. It follows that
	\begin{align*}
		|E|	&\geq	\binom{n - k}{2} + \sum_{i = 1}^{k} (d - i + 1)	\\
			&=		\binom{n - k}{2} + k(d + 1) - \frac{k(k + 1)}{2} \\
			&=		\frac{1}{2} ((n - k)(n - k - 1) + 2kd + 2k - k^{2} - k)	\\
			&=		\frac{n^{2} - 2nk + k^{2} - n + k + 2kd - k^{2} + k}{2}	\\
			&=		\frac{n^{2} - n}{2} + kd + k - nk	\\
			&=		\binom{n}{2} - k(n - d - 1).	\qedhere
	\end{align*}
\end{proof}
\begin{cor}	\label{cor:boundsOnkAndd}
	Let $\Gamma$ be an $(n, k, d)$-MBA graph. Then $n < 2k$ and ${d \leq n - \frac{k}{2k - n}}$.
\end{cor}
\begin{proof}
	We transform the inequality in \eqref{eq:boundsEdgesMBA}:
	\begin{alignat}{2}
		&\quad&	\binom{n}{2} - k(n - d - 1)	&	\leq	\frac{n(d - 1) + k}{2}	\nonumber \\
		\iff &&	n^{2} - n - 2k(n - d - 1)		&	\leq nd - n + k \nonumber	\\
		\iff &&	n^{2} - 2kn + 2kd + 2k		&	\leq nd + k		\label{eq:intermediateInequalityBoundd}\\
		\iff &&	n^{2} - nd 				&	\leq k (2n - 2d - 1). \nonumber
	\end{alignat}
	By dividing both sides of the last inequality by $2n - 2d - 1 \ne 0$, we have
	\[
		\frac{n}{2} < \frac{n(n - d)}{2(n - d) - 1} \leq k,
	\]
	where the strict inequality is equivalent with $2n(n - d) - n < 2n(n - d)$.
	
	Continuing from \eqref{eq:intermediateInequalityBoundd}, we find
	\[
		d(2k - n) \leq 2kn - k - n^{2}
	\]
	and hence
	\[
		d \leq \frac{2kn - n^{2} - k}{2k - n} = n - \frac{k}{2k - n}. \qedhere
	\]
\end{proof}
Note that since $k$ is an integer, the condition $2k > n$ is equivalent to $k \geq \ceil{\frac{n + 1}{2}}$.

\begin{lemma}	\label{lem:degreePlusVmaxAtLeastVPlus1}
	Let $\Gamma$ be an $(n, k, d)$-MBA graph. Then $k + d \geq n + 1$.
\end{lemma}
\begin{proof}
	Since $\Gamma(V \setminus \Vmax)$ is isomorphic to $K^{n - k}$, the degree of a vertex in $V \setminus \Vmax$ is at least $n - k - 1$. Since $\Gamma$ is connected, there is a $v \in V \setminus \Vmax$ adjacent to a vertex $w \in \Vmax$. For said vertex $v$, we have $d > \deg(v) \geq n - k$. Hence, $n - k \leq d - 1$, and consequently $n + 1 \leq k + d$.
\end{proof}
All results combined give some non-existence results of certain max-by-abelian graphs.
\begin{cor}	\label{cor:MBAAtLeast5Vertices}
	Let $\Gamma$ be an $(n, k, d)$-MBA graph. Then $n \geq 5$.
\end{cor}
\begin{proof}
	Combining $k \leq n - 1$ with \cref{lem:degreePlusVmaxAtLeastVPlus1}, we obtain 
	\[
		2 = n + 1 - (n - 1) \leq n + 1 - k \leq d \leq n - 2,
	\]
	hence $n \geq 4$. If $n = 4$, then $d = 2$, and $\ceil{\frac{n + 1}{2}} \leq k \leq n - 1$ implies $k = 3$. The sole vertex $v$ of $\Gamma$ not in $\Vmax$ then has degree $1$, since $\Gamma$ is connected. The sum of the degrees is then $k \cdot d + 1 = 7$. This, however, contradicts the Handshaking Lemma. We conclude that $n \geq 5$.
\end{proof}

\begin{prop}	\label{prop:MBAMinimalkImpliesNeven}
	Let $\Gamma$ be an $\left(n, \ceil{\frac{n + 1}{2}}, d\right)$-MBA graph. Then $n$ is even.
\end{prop}
\begin{proof}
	Suppose $n$ is odd. Then $k = \ceil{\frac{n + 1}{2}} = \frac{n + 1}{2}$. Applying \cref{cor:boundsOnkAndd} gives
	\[
		d \leq n - \frac{(n + 1) / 2}{2(n + 1) / 2 - n} = n - \frac{n + 1}{2} = \frac{n - 1}{2}
	\]
	and thus $k + d \leq n$, whereas \cref{lem:degreePlusVmaxAtLeastVPlus1} states that $k + d \geq n + 1$, a contradiction. Hence, $n$ is even.
\end{proof}

\subsection{\texorpdfstring{$\Rinf$-property for large values of $k$}{}}
Next, we focus on the associated RAAG of a max-by-abelian graph $\Gamma$. For large values of $k$, \ie $k$ close to the upper bound of $n - 1$, the graph is structured enough to find characteristic vertex-subgroups resulting in a non-trivial and non-abelian quotient.

We start with $(n, n - 1, d)$-MBA graphs.

\begin{prop}	\label{prop:intersectionLinkNonMaximalDegreeCharacteristic}
	Let $\Gamma$ be a graph. Put
	\[
		\tilde{V} = \bigcap_{v \notin \Vmax} lk(v).
	\]
	Then $N(\tilde{V})$ is characteristic in $A_{\Gamma}$.
\end{prop}
\begin{proof}
	If $\tilde{V}$ is empty, the claim is trivial. Since $v \notin lk(v)$ for any $v \in V$, we have that $(V \setminus \Vmax) \cap \tilde{V} = \emptyset$, hence $\tilde{V} \subseteq \Vmax$.
	
	Let $v \in \tilde{V}$. By \cref{theo:classificationCharacteristicVertexSubgroups}, it is sufficient to prove that $\Vchar(v) \subseteq \tilde{V}$. Put $V_{0} = \{v\}$ and let $V_{i}, V_{\omega}$ be as in the construction of $\Vchar(v)$, \ie
	\[
		V_{i} = \{w \in V \mid \exists v' \in V_{i - 1}: v' \leq w\}
	\]
	and $V_{\omega} = \bigcup_{n \in \N} V_{n}$. We prove by induction that $V_{i} \subseteq \tilde{V}$. For $i = 0$, this is trivial. Suppose that $V_{i} \subseteq \tilde{V}$ and that $w \in V$ dominates $v' \in V_{i}$. As $v' \in \tilde{V} \subseteq \Vmax$ by the induction hypothesis, also $w \in \Vmax$ and \(v'' \in lk(v') \subseteq st(w)\) for all $v'' \in V \setminus \Vmax$.

	Since $w \in \Vmax$, this implies that $w \in \bigcap_{v'' \notin \Vmax} lk(v'') = \tilde{V}$. Consequently, $V_{i + 1} \subseteq \tilde{V}$. We conclude that $V_{\omega} \subseteq \tilde{V}$.
	
	Now, for $\phi \in \Aut(\Gamma)$, note that $\phi(\Vmax) = \Vmax$, hence $\phi(V \setminus \Vmax) = V \setminus \Vmax$ and thus $\phi(\tilde{V}) = \tilde{V}$. In particular is
	\(
		\phi(V_{\omega}) \subseteq \phi(\tilde{V}) = \tilde{V},
	\)
	therefore, $\Vchar(v) \subseteq \tilde{V}$.
\end{proof}

\begin{prop}	\label{prop:nn-1dMBAHaveRinf}
	Let $\Gamma$ be an $(n, n - 1, d)$-MBA graph. Then $A_{\Gamma} \in \Rinf$.
\end{prop}
\begin{proof}
	Let $v$ be the sole vertex in $V \setminus \Vmax$. By \cref{prop:intersectionLinkNonMaximalDegreeCharacteristic}, $N(lk(v))$ is characteristic in $A_{\Gamma}$. Note that $lk(v) \ne \Vmax$, since
	\[
		|lk(v)| \leq d - 1 \leq n - 3 < n - 1 = |\Vmax|.
	\]
	Hence, $\tilde{\Gamma} := \Gamma(V \setminus lk(v))$ is equal to $\Gamma(\{v\}) \sqcup \Gamma(\Vmax \setminus lk(v))$. \cref{theo:RinfDisconnectedGraphs} implies that $A_{\tilde{\Gamma}}$ has the $\Rinf$-property. Since $A_{\tilde{\Gamma}} \cong \frac{A_{\Gamma}}{N(lk(v))}$ (by \cref{prop:eliminatingGenerators}) is a characteristic quotient of $A_{\Gamma}$, also $A_{\Gamma}$ has the $\Rinf$-property.
\end{proof}

Now we move on to $(n, n - 2, d)$-MBA graphs.
\begin{prop}	\label{prop:intersectionVmaxUnionLinksMBA}
	Let $\Gamma$ be a graph. Put
	\[
		\tilde{V} = \Vmax \cap \left(\bigcup_{v \notin \Vmax} lk(v)\right).
	\]
	Then $N(\tilde{V})$ is characteristic in $A_{\Gamma}$.
\end{prop}
\begin{proof}
	Let $v \in \tilde{V}$. Again, by \cref{theo:classificationCharacteristicVertexSubgroups}, it is sufficient to prove that $\Vchar(v) \subseteq \tilde{V}$. Put $V_{0} = \{v\}$ and let $V_{i}, V_{\omega}$ be as usual. We prove by induction that $V_{i} \subseteq \tilde{V}$. For $i = 0$, this is trivial. Suppose $V_{i} \subseteq \tilde{V}$ and let $w \in V$ be a vertex dominating $v' \in V_{i}$. By the induction hypothesis, $v' \in \Vmax$. Since $w$ dominates $v'$, $\deg(w) \geq \deg(v') = \Delta(\Gamma)$, so $w \in \Vmax$.
	
	Also by the induction hypothesis, there exists a $w' \notin \Vmax$ such that $v' \in lk(w')$. Consequently is $w' \in lk(v') \subseteq st(w)$. As $w' \notin \Vmax$, we have that $w' \ne w$, so $w' \in lk(w)$. Therefore, $w \in lk(w')$ implying that  $w \in \tilde{V}$. We conclude that $V_{\omega} \subseteq \tilde{V}$.
	
	Next, let $\phi \in \Aut(\Gamma)$. As $\phi(\Vmax) = \Vmax$, we have that $\phi(V \setminus \Vmax) = V \setminus \Vmax$ and hence
	\[
		\phi(\tilde{V}) = \Vmax \cap \left(\bigcup_{\phi(v) \notin \Vmax} lk(\phi(v))\right) = \Vmax \cap \left(\bigcup_{v \notin \Vmax} lk(v)\right) = \tilde{V}.
	\]
	Consequently, $\phi(V_{\omega}) \subseteq \tilde{V}$ and thus $\Vchar(v) \subseteq \tilde{V}$.
\end{proof}
Before heading to the next result, we state three commutator identities, each of which can be easily checked by working out both sides of the equality. 
\begin{lemma}	\label{lem:commutatorIdentities}
	Let $G$ be a group, and $a, b, c \in G$. Then the following identities hold:
	\begin{enumerate}[(i)]
		\item \([ab, c] = [a, c]^{b}[b, c]\)
		\item $[a^b, c] = [c, b]^{a^b} [a, c]^b [b, c]$
		\item $[a^b, c] = [b, a] [a, c] [a, b]^c$
	\end{enumerate}
\end{lemma}

\begin{prop}	\label{prop:simplifyingn-2MBAGraphs}
	Let $\Gamma$ be an $(n, n - 2, d)$-MBA graph. Write $\{v_{1}, v_{2}\} = V \setminus \Vmax$. Suppose that $V$ is the disjoint union of $lk(v_{1})$ and $lk(v_{2})$. Then
	\begin{enumerate}[(i)]
		\item $\Gamma(lk(v_{i}))$ is disconnected for $i = 1, 2$.
		\item The normal subgroup
		\(
			N := \normcl{\{[v, w] \mid v \in lk(v_{1}), w \in lk(v_{2})\}}_{A_{\Gamma}}
		\)
		is characteristic in $A_{\Gamma}$.
		\item With $N$ as above and $\tilde{\Gamma} := \Gamma(lk(v_{1})) * \Gamma(lk(v_{2}))$, we have
		\[
			\frac{A_{\Gamma}}{N} \cong A_{\tilde{\Gamma}}.
		\] 
	\end{enumerate}
\end{prop}
\begin{proof}
	We give the proof for $i = 1$, the case $i = 2$ is analogous. Since $v_{2} \notin lk(v_{2})$ and $lk(v_{1}) \cup lk(v_{2}) = V$, we have $v_{2} \in lk(v_{1})$. Let $w \in lk(v_{1}) \setminus \{v_{2}\}$ be a vertex. Because $lk(v_{1})$ and $lk(v_{2})$ are disjoint, $w$ does not lie in $lk(v_{2})$. This means that $v_{2}$ is not connected to any other vertex in $lk(v_{1})$, implying that $v_{2}$ is an isolated vertex in $\Gamma(lk(v_{1}))$. Hence, $\Gamma(lk(v_{1}))$ is disconnected.
	
	\medskip
	Let $\phi \in \Aut(A_{\Gamma})$ be an automorphism of basic type. If $\phi = \imath_{a}$ is an inversion, then for vertices $b \ne c$ different from $a$ we have
	\[
		\imath_{a}([a, b]) = [\inv{a}, b] = [b, a]^{\inv{a}}
	\]
	and $\imath_{a}([b, c]) = [b, c]$.
	Hence, if $[a, b] \in N$, then $\imath_{a}([a, b]) \in N$ as well.
	
	Next, suppose that $\phi$ is a graph automorphism. As $\phi(\Vmax) = \Vmax$, we have $\phi(\{v_{1}, v_{2}\}) = \{v_{1}, v_{2}\}$. So either $\phi(v_{i}) = v_{i}$ or $\phi(v_{i}) = v_{3 - i}$ for $i = 1, 2$. Let $v \in lk(v_{1})$ and $w \in lk(v_{2})$. In the first case, $\phi(v) \in lk(v_{1})$ and $\phi(w) \in lk(v_{2})$, so $\phi([v, w]) = [\phi(v), \phi(w)] \in N$. In the second case, $\phi(v) \in lk(v_{2})$ and $\phi(w) \in lk(v_{1})$, and hence also $\phi([v, w]) = [\phi(v), \phi(w)] \in N$.
	
	Suppose $\phi = \tau_{ab}$ is a transvection with $a \in lk(v_{1})$; the case $a \in lk(v_{2})$ is analogous and since $lk(v_{1}) \cup lk(v_{2}) = V$, this covers all cases. Let $w \in lk(v_{2})$. We have to show that $\tau_{ab}([a, w]) = [ab, w]$ lies in $N$. As $a \in lk(v_{1})$, we have $v_{1} \in lk(a) \subseteq st(b)$, where the latter inclusion follows from $a \leq b$.
	
	First, suppose that $b \ne v_{1}$. Then $v_{1} \in lk(b)$, so $b \in lk(v_{1})$, implying that $[b, w] \in N$. Thus,
	\[
		\tau_{ab}([a, w]) = [ab, w] = [a, w]^{b}[b, w] \in N.
	\]
	
	Now, if $b = v_{1}$, then $\deg(a) \leq \deg(v_{1}) < \Delta(\Gamma)$, implying that $a = v_{2}$, as $V \setminus \Vmax = \{v_{1}, v_{2}\}$. Hence, $lk(v_{2}) \subseteq st(v_{1})$. As $lk(v_{1})$ and $lk(v_{2})$ are disjoint, this implies that $lk(v_{2}) = \{v_{1}\}$. Consequently, $lk(v_{1}) = V \setminus \{v_{1}\}$ and thus $\deg(v_{1}) = n - 1$, which contradicts $v_{1} \notin \Vmax$.
	
	Finally, suppose $\phi = \gamma_{b, C}$ is a partial conjugation. Let $v \in lk(v_{1})$ and $w \in lk(v_{2})$. If $v$ and $w$ are both (not) conjugated by $b$ via $\phi$, we are done. So, suppose $\phi(v) = v^{b}$ and $\phi(w) = w$. By \cref{lem:commutatorIdentities}, we have
	\begin{align}
		[v^{b}, w]	&=	[w, b]^{v^{b}} [v, w]^{b} [b, w]		\label{eq:partialConjugationMBAn-2Case1}\\
				&=	[b, v][v, w][v, b]^{w}				\label{eq:partialConjugationMBAn-2Case2}.
	\end{align}
	If $b \in lk(v_{1})$, then \eqref{eq:partialConjugationMBAn-2Case1} implies that $[v^{b}, w] \in N$. If $b \in lk(v_{2})$, then the result follows from \eqref{eq:partialConjugationMBAn-2Case2}. This finishes the proof that $N$ is characteristic in $A_{\Gamma}$.
	
	The isomorphism is a direct application of \cref{prop:addingEdges}.
\end{proof}
\begin{cor}
	Let $\Gamma$ be an $(n, n - 2, d)$-MBA graph satisfying the same conditions as in \cref{prop:simplifyingn-2MBAGraphs}. Then $A_{\Gamma} \in \Rinf$.
\end{cor}
\begin{proof}
	We use the same notation as before. Since $\Gamma(lk(v_{1}))$ is disconnected, \cref{theo:RinfDisconnectedGraphs} implies that $A_{\Gamma(lk(v_{1}))} \in \Rinf$. Putting $\tilde{\Gamma} = \Gamma(lk(v_{1})) * \Gamma(lk(v_{2}))$, we note that this is the maximal decomposition of $\tilde{\Gamma}$, since $\Gamma(lk(v_{1}))$ and $\Gamma(lk(v_{2}))$ are both disconnected graphs. Therefore, \cref{theo:directProductsRAAGsRinf} implies that $A_{\tilde{\Gamma}} \in \Rinf$ as well. Finally, since $A_{\tilde{\Gamma}}$ is a characteristic quotient of $A_{\Gamma}$, we conclude that $A_{\Gamma}$ has the $\Rinf$-property.
\end{proof}

\begin{example}
	We provide two examples of $(n, n - 2, d)$-MBA graphs: one satisfying the conditions of \cref{prop:simplifyingn-2MBAGraphs} and one not.
	Consider the graph ${\Gamma_{1} := (K^{1} \sqcup K^{2}) * (K^{1} \sqcup K^{2})}$ (see \Cref{fig:MBAgraph6vertices}; one copy of $K^{1} \sqcup K^{2}$ has solid vertices, the other one hollow). From \cref{ex:simplicialJoinDisjointUnionCompleteGraphs}, we know that $\Gamma_{1}$ is a $(6, 4, 4)$-MBA graph.
	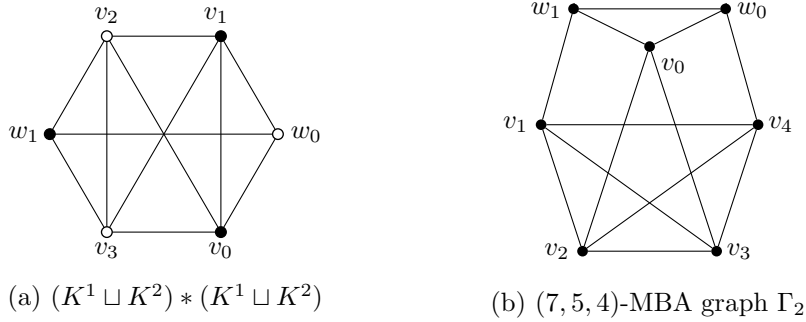
\begin{figure}[h]
		\centering
		\begin{subfigure}[h]{0.4\textwidth}
		\centering
			\begin{tikzpicture}
				\GraphInit[vstyle = Simple]
				\SetVertexMath
				\SetUpVertex[LineWidth = 0.5pt]
				\tikzset{VertexStyle/.append style = {shape = circle,fill = white,minimum size = 4pt,inner sep = 0pt}}
		        		\SetUpEdge[lw  = 0.4pt]  
				\grEmptyCycle[RA = 1.5, prefix = a]{3}
				\tikzset{VertexStyle/.append style = {fill = black}}
				\grEmptyCycle[RA = 1.5, prefix = b, rotation = 180]{3}
				\EdgeFromOneToAll{a}{b}{0}{3}
				\EdgeFromOneToAll{a}{b}{1}{3}
				\EdgeFromOneToAll{a}{b}{2}{3}
				\EdgeFromOneToSel{a}{a}{1}{2}
				\EdgeFromOneToSel{b}{b}{1}{2}
				\tikzset{AssignStyle/.append style = {right = 1pt}}
				\AssignVertexLabel[color = black, size = \small]{a}{$w_{0}$, , }
				\tikzset{AssignStyle/.append style = {left = 1pt}}
				\AssignVertexLabel[color = black, size = \small]{b}{$w_{1}$, ,}
				\tikzset{AssignStyle/.append style = {above = 1pt}}
				\AssignVertexLabel[color = black, size = \small]{a}{, $v_{2}$, }
				\AssignVertexLabel[color = black, size = \small]{b}{, ,$v_{1}$ }
				\tikzset{AssignStyle/.append style = {below = 1pt}}
				\AssignVertexLabel[color = black, size = \small]{b}{, $v_{0}$, }
				\AssignVertexLabel[color = black, size = \small]{a}{, , $v_{3}$}
			\end{tikzpicture}
			\caption{$(K^{1} \sqcup K^{2}) * (K^{1} \sqcup K^{2})$}
			\label{fig:MBAgraph6vertices}
		\end{subfigure}
		\begin{subfigure}[h]{0.4\textwidth}
		\centering
			\begin{tikzpicture}
				\GraphInit[vstyle = Simple]
				\SetVertexMath
				\tikzset{VertexStyle/.style = {shape = circle,fill = black,minimum size = 4pt,inner sep=0pt}}
			        	\SetUpEdge[lw  = 0.4pt]  
				\grComplete[y = 2, RA = 1, prefix = w, rotation = 0]{2}
				\grEmptyCycle[RA = 1.5, prefix=v, rotation = 90]{5}%
				\EdgeFromOneToSel{w}{v}{0}{0,4}
				\EdgeFromOneToSel{w}{v}{1}{0,1}
				\tikzset{AssignStyle/.append style = {left = 1pt}}
				\AssignVertexLabel[color = black, size = \small]{v}{, $v_{1}$, $v_{2}$, ,}
				\tikzset{AssignStyle/.append style = {right = 1pt}}
				\AssignVertexLabel[color = black, size = \small]{v}{, , , $v_{3}$ ,$v_{4}$}
				\tikzset{AssignStyle/.append style = {below right = .5pt}}
				\AssignVertexLabel[color = black, size = \small]{v}{$v_{0}$, , , ,}
				\tikzset{AssignStyle/.append style = {left = 1pt}}
				\AssignVertexLabel[color = black, size = \small]{w}{, $w_{1}$}
				\tikzset{AssignStyle/.append style = {right = 1pt}}
				\AssignVertexLabel[color = black, size = \small]{w}{$w_{0}$, }
				\EdgeInGraphSeq{v}{1}{3}
				\EdgeInGraphMod{v}{5}{2}
			\end{tikzpicture}
			\caption{$(7, 5, 4)$-MBA graph $\Gamma_{2}$}
			\label{fig:MBAgraph7vertices}
		\end{subfigure}
		\caption{Further examples of MBA-graphs}
	\end{figure}
	
	The vertices in $V \setminus \Vmax$ are $w_{0}$ and $w_{1}$, for which we have
	\[
		lk(w_{0}) = \{w_{1}, v_{0}, v_{1}\}, \quad lk(w_{1}) = \{w_{0}, v_{2}, v_{3}\},
	\]
	showing that $\Gamma_{1}$ satisfies the conditions of \cref{prop:simplifyingn-2MBAGraphs}.
	
	On the other hand, consider $\Gamma_{2}$ as in \Cref{fig:MBAgraph7vertices}. Then $V \setminus \Vmax = \{w_{0}, w_{1}\}$, but $v_{0}$ is adjacent to both $w_{0}$ and $w_{1}$. Moreover, $v_{2}$ lies in neither $lk(w_{0})$ nor $lk(w_{1})$.
\end{example}

For $(n, n - 2, d)$-MBA graphs not satisfying the conditions of \cref{prop:simplifyingn-2MBAGraphs}, we can nonetheless find a suitable characteristic vertex subgroup.

\begin{lemma}	\label{lem:MBAIntersectionsUnionsLinks}
	Let $\Gamma$ be a max-by-abelian graph. Then
	\begin{equation}	\label{eq:MBAintersectionLinksSubsetVmax}
		\bigcap_{v \notin \Vmax} lk(v) \subsetneq \Vmax
	\end{equation}
	and
	\begin{equation}	\label{eq:MBAunionIntersectsVmax}
		\Vmax \cap \left(\bigcup_{v \notin \Vmax} lk(v)\right) \ne \emptyset.
	\end{equation}
\end{lemma}
\begin{proof}
	In \cref{prop:intersectionLinkNonMaximalDegreeCharacteristic}, we have proven that the inclusion of \eqref{eq:MBAintersectionLinksSubsetVmax} holds, albeit not necessarily strict. Suppose equality holds. Then for ${v \in V \setminus \Vmax}$, we have $\Vmax \subseteq lk(v)$. As $\Gamma(V \setminus \Vmax)$ is complete, we also have
	\[
		V \setminus (\Vmax \cup \{v\}) \subseteq lk(v).
	\]
	Consequently, $V \setminus \{v\} \subseteq lk(v)$, contradicting $v \notin \Vmax$. Hence, the inclusion in \eqref{eq:MBAintersectionLinksSubsetVmax} is strict.
	
	For \eqref{eq:MBAunionIntersectsVmax}, note that since $\Gamma$ is connected, there is a $v \notin \Vmax$ that is adjacent to a vertex in $\Vmax$. Hence, the intersection in \eqref{eq:MBAunionIntersectsVmax} is non-empty.
\end{proof}
\begin{prop}	\label{prop:characteristicQuotientn-2MBAGraph}
	Let $\Gamma$ be an $(n, n - 2, d)$-MBA graph and put $V \setminus \Vmax = \{v_{1}, v_{2}\}$.
	\begin{enumerate}[(i)]
		\item If $lk(v_{1}) \cap lk(v_{2}) \ne \emptyset$, then $\Gamma\big(V \setminus (lk(v_{1}) \cap lk(v_{2}))\big)$ is non-complete.	\label{item:nonDisjointLinks}
		\item If $lk(v_{1}) \cup lk(v_{2}) \ne V$, then $\Gamma\big(V \setminus (\Vmax \cap (lk(v_{1}) \cup lk(v_{2})))\big)$ is non-complete.	\label{item:nonCoveringLinks}
	\end{enumerate}
	In either case, $A_{\Gamma}$ has a non-abelian RAAG $A_{\Gamma'}$ as characteristic quotient with $|V'| < |V|$.
\end{prop}
\begin{proof}
	Suppose $lk(v_{1}) \cap lk(v_{2}) \ne \emptyset$. Note that $v_{1}, v_{2} \in V \setminus (lk(v_{1}) \cap lk(v_{2}))$ and because $lk(v_{1}) \cap lk(v_{2}) \subsetneq \Vmax$ by \cref{lem:MBAIntersectionsUnionsLinks}, there is a vertex $v \in \Vmax \setminus (lk(v_{1}) \cap lk(v_{2}))$. Then by definition of $v$, either $vv_{1} \notin E$ or $vv_{2} \notin E$, showing that $\Gamma(V \setminus (lk(v_{1}) \cap lk(v_{2}))$ is non-complete.
	
	On the other hand, suppose $lk(v_{1}) \cup lk(v_{2}) \ne V$. As both $v_{1}$ and $v_{2}$ lie in $lk(v_{1}) \cup lk(v_{2})$, there is a $v \in \Vmax$ such that $v \notin lk(v_{1}) \cup lk(v_{2})$. Then $v, v_{1}$ and $v_{2}$ all lie in $V \setminus (\Vmax \cap (lk(v_{1}) \cup lk(v_{2})))$ and by definition of $v$ is $vv_{1} \notin E$. Hence, $\Gamma(V \setminus (\Vmax \cap (lk(v_{1}) \cup lk(v_{2}))))$ is non-complete.
	
	\medskip
	
	The claim regarding the characteristic quotients follows for \eqref{item:nonDisjointLinks} from \cref{lem:MBAIntersectionsUnionsLinks} and \cref{prop:intersectionLinkNonMaximalDegreeCharacteristic}. For \eqref{item:nonCoveringLinks}, it follows from \cref{lem:MBAIntersectionsUnionsLinks} and \cref{prop:intersectionVmaxUnionLinksMBA}.
\end{proof}
\section{Partial answer to conjecture}
\addtocounter{subsection}{1}

Combining all the obtained results, we are able to partially answer the conjecture.
\begin{theorem}
	Let \(\Gamma\) be a non-complete graph on at most \(7\) vertices. Then \(A_{\Gamma}\) has the \(\Rinf\)-property.
\end{theorem}
The proof boils down to applying the Simplification Lemma and showing that the results from \cref{sec:freeProductsRAAGs,sec:regularGraphs,sec:MBA} cover all remaining cases when \(\Gamma\) has at most \(7\) vertices, see \cite[p.~155]{SendenThesis} for the details.

\rmfamily

\printbibliography[heading=bibintoc]

\end{document}